\documentclass[moor]{informs1-custom}
\RequirePackage[OT1]{fontenc}
\usepackage[normalem]{ulem}
\usepackage{nicefrac,bm,upgreek,mathtools}
\usepackage[mathscr]{eucal}

\usepackage{natbib}
 \NatBibNumeric
 \def\bibfont{\small}%
 \def\BIBand{and}%
 \bibpunct[, ]{[}{]}{,}{n}{}{,}%
 \usepackage[colorlinks=true,breaklinks=true,bookmarks=true,urlcolor=blue,
     citecolor=blue,linkcolor=blue,bookmarksopen=false,draft=false]{hyperref}

\def\EMAIL#1{\href{mailto:#1}{#1}}

\TheoremsNumberedBySection
\EquationsNumberedBySection 

\usepackage[capitalize,nameinlink]{cleveref}
\newcommand{\ttup}[1]{\textup{(}#1\textup{)}}

\crefname{section}{Section}{Sections}
\crefname{subsection}{Section}{Sections}
\crefname{hypothesis}{Hypothesis}{Hypotheses}
\crefname{assumption}{Assumption}{Assumptions}
\crefname{lemma}{Lemma}{Lemmas}
\crefname{corollary}{Corollary}{Corollaries}
\crefname{theorem}{Theorem}{Theorems}
\crefname{remark}{Remark}{Claims}

\Crefname{figure}{Figure}{Figures}

\crefformat{equation}{\textup{#2(#1)#3}}
\crefrangeformat{equation}{\textup{#3(#1)#4--#5(#2)#6}}
\crefmultiformat{equation}{\textup{#2(#1)#3}}{ and \textup{#2(#1)#3}}
{, \textup{#2(#1)#3}}{, and \textup{#2(#1)#3}}
\crefrangemultiformat{equation}{\textup{#3(#1)#4--#5(#2)#6}}%
{ and \textup{#3(#1)#4--#5(#2)#6}}{, \textup{#3(#1)#4--#5(#2)#6}}%
{, and \textup{#3(#1)#4--#5(#2)#6}}

\Crefformat{equation}{#2Equation~\textup{(#1)}#3}
\Crefrangeformat{equation}{Equations~\textup{#3(#1)#4--#5(#2)#6}}
\Crefmultiformat{equation}{Equations~\textup{#2(#1)#3}}{ and \textup{#2(#1)#3}}
{, \textup{#2(#1)#3}}{, and \textup{#2(#1)#3}}
\Crefrangemultiformat{equation}{Equations~\textup{#3(#1)#4--#5(#2)#6}}%
{ and \textup{#3(#1)#4--#5(#2)#6}}{, \textup{#3(#1)#4--#5(#2)#6}}%
{, and \textup{#3(#1)#4--#5(#2)#6}}

\crefdefaultlabelformat{#2\textup{#1}#3}

\newcommand{\process}[1]{{\{#1(t)\}_{t\ge0}}}

\newcommand{\mx}{{_\mathsf{max}}}
\newcommand{\mn}{{_\mathsf{min}}}

\newcommand{\df}{\coloneqq}
\DeclareMathOperator{\Exp}{\mathbb{E}} 
\newcommand{\D}{\mathrm{d}} 
\newcommand{\E}{\mathrm{e}} 
\newcommand{\RR}{\mathbb{R}} 
\newcommand{\Rd}{{\mathbb{R}^m}} 
\newcommand{\NN}{\mathbb{N}} 
\newcommand{\ZZ}{\mathbb{Z}} 
\newcommand{\Ind}{\bm{1}} 
\newcommand{\Act}{\mathbb{U}} 
\newcommand{\Cc}{\mathcal{C}}
\newcommand{\dd}{\mathfrak{d}}

\newcommand{\Usm}{\mathfrak{U}_{\mathsf{sm}}} 
\newcommand{\transp}{^{\mathsf{T}}} 
\newcommand{\comp}{^{\mathsf{c}}} 

\newcommand{\Lyap}{{\mathscr{V}}}

\newcommand{\veo}{{\varepsilon}}
\newcommand{\cZn}{{\mathcal{Z}^n}} 
\newcommand{\fZn}{{\mathfrak{Z}^n_{\mathsf{sm}} }} 

\newcommand{\abs}[1]{\lvert#1\rvert}
\newcommand{\norm}[1]{\lVert#1\rVert}
\newcommand{\babs}[1]{\bigl\lvert#1\bigr\rvert}

\newcommand{\babss}[1]{\biggl\lvert#1\biggr\rvert}
\newcommand{\Babss}[1]{\Biggl\lvert#1\Biggr\rvert}
\newcommand{\bnorm}[1]{\bigl\lVert#1\bigr\rVert}

\DeclareMathOperator{\diag}{diag}
\DeclareMathOperator{\trace}{trace}

\newcommand{\sK}{{\mathscr{K}}} 
\newcommand{\cK}{{\mathcal{K}}} 

\newcommand{\sB}{{\mathscr{B}}} 
\newcommand{\cB}{{\mathcal{B}}} 

\newcommand{\Ag}{\mathcal{A}}
\newcommand{\sA}{\mathscr{A}}

\newcommand{\Lg}{{\mathcal{L}}}
\newcommand{\cI}{{\mathcal{I}}}
\newcommand{\sV}{{\mathscr{V}}}
\newcommand{\sX}{{\mathscr{X}}}
\newcommand{\cV}{{\mathcal{V}}}

\usepackage{accents}
\newlength{\dhatheight}

\makeatletter
\DeclareRobustCommand\widecheck[1]{{\mathpalette\@widecheck{#1}}}
\def\@widecheck#1#2{%
 \setbox\z@\hbox{\m@th$#1#2$}%
 \setbox\tw@\hbox{\m@th$#1%
 \widehat{%
 \vrule\@width\z@\@height\ht\z@
 \vrule\@height\z@\@width\wd\z@}$}%
 \dp\tw@-\ht\z@
 \@tempdima\ht\z@ \advance\@tempdima2\ht\tw@ \divide\@tempdima\thr@@
 \setbox\tw@\hbox{%
 \raise\@tempdima\hbox{\scalebox{1}[-1]{\lower\@tempdima\box
\tw@}}}%
 {\ooalign{\box\tw@ \cr \box\z@}}}
\makeatother

\begin{document}

\TITLE{\Large On uniform exponential ergodicity of Markovian multiclass
many-server queues in the Halfin--Whitt regime}
\MANUSCRIPTNO{MOR-2018-332.R1}
\RUNTITLE{Uniform exponential ergodicity of multiclass many-server queues}
\RUNAUTHOR{A. Arapostathis, H. Hmedi and G. Pang}

\ARTICLEAUTHORS{
\AUTHOR{Ari Arapostathis, Hassan Hmedi}
\AFF{Department of Electrical and Computer Engineering,
The University of Texas at Austin,\\ 2501 Speedway, EER 7.824,
Austin, TX~~78712,
\EMAIL{[ari,hmedi]@utexas.edu}}
\AUTHOR{Guodong Pang}
\AFF{The Harold and Inge Marcus Dept.
of Industrial and Manufacturing Eng.,
College of Engineering,\\
Pennsylvania State University,
University Park, PA~~16802,
\EMAIL{gup3@psu.edu}}
}

\ABSTRACT{
We study ergodic properties of Markovian
multiclass many-server queues which are uniform over scheduling policies,
as well as the size of the system. 
The system is heavily loaded in the Halfin--Whitt regime, and the scheduling
policies are work-conserving and preemptive. We provide a unified approach via 
a Lyapunov function method that establishes Foster-Lyapunov equations for both 
the limiting diffusion and the prelimit diffusion-scaled queueing
processes simultaneously.

We first study the limiting controlled diffusion, and show
that if the spare capacity (safety staffing) parameter is positive,
the diffusion is exponentially ergodic uniformly over all
stationary Markov controls, and the invariant probability measures have
 uniform exponential tails.
This result is sharp, since when there is no abandonment and the
spare capacity parameter is negative, then the controlled diffusion is transient
under any Markov control.
In addition, we show that if all the abandonment rates are positive,
the invariant probability measures have sub-Gaussian tails, regardless whether
the spare capacity parameter is positive or negative.
 
Using the above results, we proceed to establish the corresponding ergodic properties
for the diffusion-scaled queueing processes.
In addition to providing a simpler proof of the results in Gamarnik and Stolyar
[\emph{Queueing Syst} (2012) 71:25--51], we extend these results
to the multiclass models with renewal arrival processes, albeit under
the assumption that the mean residual life functions are bounded.
For the Markovian model with Poisson arrivals, we obtain stronger results
and show that the convergence to the stationary
distribution is at an exponential rate uniformly over all
work-conserving stationary Markov scheduling policies.}

\MSCCLASS{90B22, 60K25, 90B15}

\maketitle

\section{Introduction.}
Multiclass many-server queues in the Halfin--Whitt (H--W) regime have been 
extensively studied as a useful model for large-scale service systems. 
In this paper we focus on ergodic properties of such multiclass 
queueing models. 
The ergodic properties of these systems have
been the subject of great interest in applied probability
(for a discussion see
\cite{GS-12,Stolyar-15,Stolyar-Yudovina-12,Stolyar-Yudovina-13,Stolyar-15b,Dai-He-13}).
It is important to understand if a queueing system is stable 
 and the rate at which a performance measure converges to 
 the steady state under different scheduling or routing policies. 
 For the multiclass ``V" network, Gamarnik and Stolyar \cite{GS-12} prove 
 the tightness of the stationary distributions
of the diffusion-scaled state processes under any work conserving
scheduling policy, provided that there is $\sqrt{n}$ safety staffing 
($n$ is the scaling parameter).
 They show that the diffusion-scaled
queueing processes are ergodic under all
work conserving scheduling policies, and have exhibited exponential
tail bounds 
for the stationary distribution.
The proofs of these significant results utilize some natural test functions
based on the total workload,
but there is no uniform Foster--Lyapunov equation to 
exhibit the rate of convergence to the stationary distribution. 
For the limiting diffusion of the ``V" network, 
when the control equals $(0,\dots, 0,1)\transp$,
which arises as the limit of a static priority policy,
the ergodic properties established in Dieker and Gao \cite{DG-13} for a class of 
piecewise Ornstein--Uhlenbeck (O--U) processes arising in 
many-server queues with phase-type service times can be applied. 
Exponential ergodicity is also established for the limiting diffusion 
(as a special case of a more general class of SDEs)
under any constant control in Arapostathis et al. \cite{APS18}. 

\smallskip

The following important open questions are addressed in this paper: 

\begin{enumerate}
\item[(1)] 
Is the limiting controlled diffusion exponentially ergodic under
all stationary Markov controls?
How different are the tail asymptotics of the invariant measures
with or without abandonment?
\item[(2)] Is there a unified approach based on 
Foster--Lyapunov theory that can be used
to establish uniform exponential ergodicity for both the limiting diffusion
and the diffusion-scaled queueing processes?
\end{enumerate}

\smallskip

We provide affirmative answers to all these questions.
We consider multiclass models with (delayed) renewal arrivals,
class-dependent exponential service times, and 
class-dependent exponential patience times. 
We assume that the system is operating under work-conserving
and preemptive scheduling policies.
It is well known that the diffusion-scaled queueing processes under
such scheduling policies converge weakly to a limiting diffusion with a drift
given in \cref{E-drift} and a diagonal constant covariance matrix
(see \cite{AMR04,HZ-04}).

We start with the limiting controlled diffusion. 
When the controls are constant, 
the limiting diffusion has a piecewise linear drift
and belongs to a class of piecewise
O--U processes. 
Applying \cite[Theorem~3]{DG-13} to our model with positive abandonment rate,
one can deduce that the limiting diffusion is exponentially ergodic under
a specific constant control corresponding to a static priority scheduling policy
(see Remark~\ref{R2.1}).
On the other hand, it is shown in \cite[Theorem~3.5]{APS18} that the limiting diffusion
is exponentially ergodic under any constant control (see \cref{R2.2}). 
The proofs of these results rely on the construction
of a common quadratic-type Lyapunov function for the piecewise linear equations.
However, this methodology only works for constant controls, and leaves the question
of stability over Markov controls open.

We exploit Lyapunov functions that are constructed in an intricate manner
in order to capture both
the total workload on the positive half-space and the idleness on the negative
half-space.
Such functions are of course quite natural, and have been
used as test functions in \cite{GS-12} to derive tail bounds.
However, for the diffusion,
the total workload and idleness need to be treated with the proper
``weights'' or ``tilting'', 
interacting with a ``smoothing'' cut-off function
which needs to be deployed.
Such delicate care is not only needed for the drift as usual,
but more importantly, for the second-order derivatives.
For multiclass queueing models in the Halfin--Whitt regime,
such constructions appear to be necessary in order to deal with 
both the workload and idleness processes simultaneously. 
This constitutes our first main methodology contribution in this paper. 

We present Foster--Lyapunov equations that are uniform over
all Markov controls, and show that 
\begin{enumerate}
\item[(a)] if the spare capacity parameter (safety staffing)
is positive, then the limiting diffusion is uniformly exponentially ergodic,
and the corresponding invariant probability measures have
 uniform exponential tails;
\item[(b)]
when the abandonment rates are all positive, regardless the spare capacity
parameter being positive or negative, in addition to uniform exponential
ergodicity, we show that the invariant probability measures have sub-Gaussian tails.
\end{enumerate}
These answer the questions in (1) above. 

We then show that the Foster--Lyapunov equations for the limiting diffusion
offer a very natural tool with which we establish
uniform ergodic properties for the diffusion-scaled queueing processes. 
This answers the question in (2). 
In this manner we provide a \emph{unified} approach to  the study
of the limiting diffusion and the corresponding diffusion-scaled processes. 

In the case of Poisson arrivals, by employing the same Lyapunov functions
used for the limiting diffusion, we show that the corresponding results
in (a) above hold for the diffusion-scaled queueing processes
(see also \cref{S2.3}).
Sub-Gaussian tails are not possible for the invariant distribution
of the diffusion-scaled queueing processes, and one can only hope for tails
that decay faster than any exponential. On the other hand,
when the abandonment rates are all positive, we improve somewhat on the
results in \cite{GS-12}, although a conjecture stated in that paper
still remains open. 
Even though in the cases of Poisson and renewal arrivals the limiting diffusions
agree,
with the only differences lying in the covariance functions, 
for the analysis of the prelimit processes, 
we need to augment the state process in the renewal case. 

With renewal arrivals, we consider the Markov process composed of the diffusion-scaled 
queueing processes and the interarrival age processes of the renewal arrivals.
The Lyapunov functions used for the limiting diffusion are adapted to construct
 appropriate Lyapunov functions for the joint processes.
On the other hand, the hazard rate functions and mean residual lifetime functions
of the interarrival times must be also used in a proper manner to take into account 
the age processes as suggested in \cite{Takis-99}. 
We prove the following results under the assumption that
the residual lifetime function is bounded:
(a$^\prime$) if the spare capacity parameter (safety staffing)
is positive, we prove a Foster--Lyapunov equation, which shows that the
joint Markov process is positive Harris recurrent under any work-conserving
stationary Markov scheduling policy; 
(b$^\prime$) if the abandonment rates are all positive,  we obtain
a Foster--Lyapunov equation which shows that the first absolute moments
of the invariant distribution are uniformly bounded.
If we impose the additional assumption that the hazard rate function is bounded,
we show that the marginal of the stationary distribution corresponding
to the queueing state has exponential moments. 

This work also relates to the vast literature on the validity of diffusion
approximations for queues in heavy traffic. 
We focus on the literature of many-server queueing models in the
H--W regime, and refer the readers to
\cite{GZ-06, Budhiraja-09, Katsuda-10, Gurvich-16, YY-16, YY-18, BDM-17}
and references therein for results in the conventional (single-server)
heavy-traffic regime. 
For the single-class $GI/M/n$ queues, \citet{Halfin-Whitt}
established the interchange of limits, and they used a bounding argument via
single-server queues to show the tightness of the steady-state distributions
of the diffusion-scaled processes. 
\citet{Dai-14} studied the validity of the multidimensional
diffusion approximations for $GI/Ph/n+M$ queues with phase-type service times. 
\citet{AR-17} proved the convergence of the stationary
distributions of suitably scaled infinite-dimensional measure-valued
processes for the $GI/GI/N$ queues in the H--W regime,
and they also studied the ergodic properties of the SPDE limit of the
same model in \cite{AR-15}. We also refer the readers to the steady state
analysis of many-server queues in \cite{GG-13, GM-08, BD-17, BDF-17}. 
All these studies are on the single-class many-server queues. 
For multiclass many-server queues in the
H--W regime this topic still remains wide open. 
The only known result is for the Markovian `N' network
\cite{Stolyar-15} where Stolyar proves the
interchange of limits for the model without abandonment
under a particular static priority policy. 

Uniform exponential ergodicity can substantially simplify
the study of ergodic control problems, since there is a rich body of
existing theory that can be applied \cite[Chapter~3]{book}.
On the other hand, if the system is not endowed with such blanket stability properties,
and the running cost functional is not near-monotone,
then the analysis of these problems can be quite involved.
In the study of ergodic control of the ``V" network in \cite{ABP15}, 
a key structural property of the system dynamics had to be identified 
due to the lack of uniform stability and near-monotonicity of the running cost. 
It was assumed that all the abandonment rates are strictly positive but
no positive safety staffing requirement was imposed.
The results in this paper enable the study of ergodic control problems
for the ``V" network when there is no abandonment
but there is positive safety staffing. 
Uniform stability properties are yet to be explored for multiclass multi-pool networks. 
Without such blanket stability properties, ergodic control problems for
multiclass multi-pool networks have been recently studied in
\cite{AP-16, AP-18-MOR, AP-18-SPA}, under the hypothesis that at least one
abandonment rate is positive. 

\subsection{Notation.}

We summarize some of the notation used throughout the paper. 
We use $\Rd$ (and $\mathbb{R}^m_+$), $m\ge 1$, to denote real-valued
$m$-dimensional (nonnegative) vectors, and write $\RR$ for $m=1$.
For $x, y\in \RR$, we let
$x \vee y = \max\{x,y\}$, $x\wedge y = \min\{x,y\}$, 
$x^+ = \max\{x, 0\}$ and $x^- = \max\{-x,0\}$. 
For a set $A\subseteq\Rd$, we use
 $A^{\mathsf c}$, $\partial A$, and $\Ind_{A}$ to denote 
the complement, the boundary, and the indicator function of $A$, respectively.
A ball of radius $r>0$ in $\Rd$ around a point $x$ is denoted by $\sB_{r}(x)$,
or simply as $\sB_{r}$ if $x=0$.
We also let $\sB \equiv \sB_{1}$.
The Euclidean norm on $\Rd$ is denoted by $\abs{\,\cdot\,}$,
and $\langle \cdot\,,\,\cdot\rangle$ stands for the inner product.
Also for $x\in\Rd$, we let $\norm{x}^{}_1\df \sum_i \abs{x_i}$,
$x_\mx\df \max_i\, x_i$, and $x_{\mathsf{min}}\df \min_i\, x_i$,
and $x^\pm\df \bigl(x_1^\pm,\dotsc,x_m^\pm\bigr)$.
For a finite signed measure $\nu$ on $\Rd$,
and a Borel measurable $f\colon\Rd\to[1,\infty)$,
we define the $f$-norm of $\nu$ by
\begin{equation}\label{E-fnorm}
\norm{\nu}_f \,\df\, \sup_{\substack{g\in\cB(\Rd), \; \abs{g}\le f}}\;
\babss{\int_{\Rd} g(x)\,\nu(\D{x})}\,,
\end{equation}
where $\cB(\Rd)$ denotes the class of Borel measurable functions on $\Rd$.
Observe that $\norm{\,\cdot\,}_1=\norm{\,\cdot\,}_{\mathsf{TV}}$, the latter denoting
the total variation norm.

\section{Uniform exponential ergodicity of the diffusion limit.}

In \cref{S2.1} we describe the limiting diffusion,
and proceed with a summary of the results and the technical approach
in \cref{S2.2,S2.3}, respectively.
Some important definitions are in \cref{S2.4}, followed by the
main technical results in \cref{S2.5,S2.6}.

\subsection{The limiting controlled diffusion.}\label{S2.1}

We consider a controlled
$m$-dimensional stochastic differential equation (SDE) of the form
\begin{equation}\label{E-sde}
\D X_t\,=\,b(X_t,U_t)\,\D{t} +\upsigma(X(t))\, \D W_t\,,\qquad X(0)=x_0\in\Rd\,,
\end{equation} 
with
$b\colon\Rd\to\Rd$ given by
\begin{equation}\label{E-drift}
b(x,u) \,=\, \ell-M(x-\langle  e,x\rangle^+u)-\langle  e,x\rangle^+\varGamma u
\,=\, \begin{cases}
\ell - \bigl(M +(\varGamma-M)u e\transp\bigr) x\,, & \langle e,x\rangle>0\,,\\[2pt]
\ell - Mx\,, &\langle e,x\rangle\le0\,.
\end{cases}
\end{equation}
Here, $\ell \in \Rd$, $u\in\RR^m_+$ satisfies $\langle  e,u\rangle=1$
with $e = (1,\dotsc,1)\transp\in\Rd$, $M=\diag(\mu_1,\dotsc,\mu_m)$ is a
positive diagonal matrix,
and $\varGamma=\diag(\gamma_1,\dotsc,\gamma_m)$ with $\gamma_i\in\RR_{+}$,
$i\in \mathcal{I}\df\{1,\dotsc,m\}$.
The process $W_t$ is a standard $m$-dimensional Brownian motion, and
the covariance function $\upsigma\colon \Rd \to \RR^{m\times m}$ is a
positive diagonal matrix.
Such a process arises as a limit of the suitably scaled queueing processes of
multiclass Markovian many-server queues in the H--W regime \cite{AMR04, HZ-04}. 

In these models, if the scheduling policy is based on a static priority assignment
on the queues, then the vector $u$ in \cref{E-sde} corresponds to a constant control
which is
an extreme point of the convex set
\begin{equation*}
\varDelta \,\df\, \{u\in\RR^m \,\colon u\ge0\,,\ \langle e,u\rangle = 1\}\,.
\end{equation*}

\begin{remark} \label{R2.1}
As mentioned earlier, ergodicity and exponential ergodicity
of a class of piecewise O--U processes as in \cref{E-sde}
have been addressed in \cite{DG-13}. 
In this model, they assume that $M$ is a nonsingular M-matrix
such that the vector $e\transp M$ has nonnegative components,
$\varGamma=\alpha I$, and $\ell=-\beta u$ for positive constants $\alpha$, $\beta$, 
and a constant vector $u \in \varDelta$.
Applying their results to the multiclass $M/M/n+M$ model with abandonment, 
exponential ergodicity of the limiting diffusion under the specific constant control
$\bar{u} = (0,\dots,0,1)\transp$, corresponding to class $m$ being given the
least priority, 
is established in \cite[Theorem 3]{DG-13}. 
On the other hand, for the multiclass $M/M/n$ model without abandonment, that is,
$\varGamma=0$, positive recurrence is established for the limiting
diffusion with the control $\bar{u}$
but the rate of convergence is not identified \cite[Theorem 2]{DG-13}. 
\end{remark}

\begin{remark} \label{R2.2}
The model in \cref{E-sde} with $M$ a nonsingular M-matrix,
and for constant control $U_t$ has also been studied extensively
in \cite{APS18} (as a special class of the L{\'e}vy--driven SDEs studied there). 
It is shown in that paper that when $\varGamma=0$, the quantity
\begin{equation}\label{E-varrho}
\varrho\,\df\, -\,\bigl\langle  e,M^{-1}\ell\,\bigr\rangle
\end{equation}
plays a fundamental role in the characterization of stability.
Specifically, it is shown in \cite[Theorem~3.2]{APS18} that
if $\varrho>0$, then $X_t$ is positive recurrent under
any constant control $U_t$, and if $\varrho<0$ ($\varrho=0$), then it is transient
(cannot be positive recurrent) under
any stationary Markov control satisfying $\varGamma v(x)=0$ a.e.
\cite[Theorem~3.3]{APS18}.
Another interesting property of $\varrho$ which we find
in \cite[Corollary~5.1]{APS18} is that,
provided $\varGamma=0$, and the diffusion
under some stationary Markov control $v$ is positive recurrent
with invariant probability measure $\uppi_v$, then necessarily
\begin{equation}\label{E-nice}
\varrho\,=\,\int_\Rd \langle e,x\rangle^-\,\uppi_v(\D{x})\,.
\end{equation}
This can be interpreted as follows:
the `average idleness' in the steady state always equals
the spare capacity parameter.
These results of course apply to the problem at hand since $M$ is
a diagonal matrix.
In addition, the rate of convergence is shown to be exponential if either
$\varGamma u=0$ or $\varGamma u\neq 0$ for any constant control $u \in \varDelta$
\cite[Corollary~4.2]{APS18}. 
\end{remark}

Let $\Usm$ denote the class of Borel measurable maps
$v\colon\Rd\to\varDelta$.
Every element $v$ of $\Usm$ is identified with
the stationary Markov controls $U_t=v(X_t)$.
Under any such control, it is well known that \cref{E-sde} has a unique strong
solution which is a strong Feller process \cite{Gyongy-96}.
Let $P^v_t(x,\D{y})$ denote its transition probability.

The diffusion in \cref{E-sde} is called \emph{uniformly stable}
(in the sense of \cite[Definition~3.3.3]{book}), if under any $v\in\Usm$,
the process $X_t$
is positive recurrent and the collection of invariant probability measures
is tight.
We say that \cref{E-sde} is \emph{uniformly exponentially ergodic}, if
it is uniformly stable and
there exist positive constants $C$ and $\gamma$ and a function
$\Lyap \colon \Rd\to[1,\infty)$ such that
\begin{equation*}
\bnorm{P^v_t(x,\cdot\,)-\uppi_v(\cdot)}_{\mathsf{TV}}\,\le\,
C \Lyap(x) \E^{-\gamma t}\quad \forall (x,t)\in\Rd\times\RR_+\,,
\end{equation*}
and all $v\in\Usm$.

\subsection{Brief summary of the results.}\label{S2.2}
In \cref{T2.1} we show that if $\varrho>0$, then 
\cref{E-sde} is uniformly exponentially ergodic.
Therefore, when $\varGamma=0$, \cref{E-nice} holds
over all stationary Markov controls $v\in\Usm$.
In addition, the invariant probability measures have
 uniform exponential tails,
and by that, we mean that there exists some $\veo>0$, such that 
$\sup_{v\in\Usm}\,\int_\Rd \E^{\veo\abs{x}}\,\uppi_v(\D{x})<\infty$.
On the other hand, if $\varGamma>0$, then the
associated invariant probability measures have sub-Gaussian tails,
that is, $\sup_{v\in\Usm}\,\int_\Rd \E^{\veo\abs{x}^2}\,\uppi(\D{x})<\infty$ for
some $\veo>0$ (see \cref{T2.2}).

In \cref{S3} we address the $n$-server networks.
We first present the results for the models with (delayed) renewal arrival processes 
in \cref{S3.2}.
The counterpart of \cref{T2.1} here is given in \cref{T3.1},
and this is established for renewal arrivals
(this should be compared to \cite[Theorem~2]{GS-12}).
In this theorem, the hazard rate functions are assumed bounded.
This is a rather strong assumption, but the result, which asserts
 uniform exponential tails for the invariant distributions
under work-conserving stationary Markov
policies is equally strong. 
With strictly positive abandonment parameters, and with the hazard rate
function only locally bounded, we establish uniform stability
of the queueing system under all work-conserving stationary Markov
policies in \cref{T3.2}. 
With possibly zero abandonment in all classes, and with positive
$\sqrt n$ safety staffing, we show in \cref{T3.3}, that the
combined renewal age and queueing state process is positive Harris recurrent.
In addition, if the limit of the safety staffing is positive,
the invariant probability distributions are tight. 
In this result, the hazard rate
function is assumed only locally bounded.

Networks with Poisson arrivals are studied in \cref{S3.3}.
We show in \cref{C3.1} that
positive spare capacity implies exponential ergodicity.
However, as noted in \cite{GS-12} the invariant distribution of
an $n$-server network cannot have a sub-Gaussian tail.
This property is recovered only at the weak limit as $n\to\infty$,
and it is worthwhile comparing \cite[Theorem~4]{GS-12} with \cref{T2.2},
which in addition shows uniform exponential ergodicity.
When all abandonment rates are positive, we can only
exhibit a stronger Foster--Lyapunov equation (see \cref{T3.4})
which implies that $\E^{\delta\abs{x}}$ is uniformly integrable over
the invariant probability distributions for any $\delta>0$.

In addition to these results, we investigate the properties in
Theorem~2\,(i) and Theorem~4\,(i) of \cite{GS-12}.
We provide proofs of the analogous results for the limiting diffusion
in \cref{L2.2,T2.3}, respectively, using Foster--Lyapunov techniques.
The counterpart of \cref{L2.2} for the $n$-system is
given in \cref{T3.5}
 and is an improvement over the statement in \cite[Theorem~2\,(i)]{GS-12}.
However, we have not been able to prove or disprove the related conjecture
in \cite[p.~33]{GS-12}.

\subsection{Summary of the technical approach.}\label{S2.3}

The first important step in the study of this problem is the
construction of appropriate Lyapunov functions.
We use two building blocks for these functions:
one represents the total workload, and the other is a measure of idleness.
The scaling of these in \cref{E-Psi} plays a crucial role.
Two scaling parameters are used: $\theta$ to balance the mix of workload
and idleness, and $\veo$ to handle the terms arising from the second derivatives
in the extended generator of the controlled diffusion.
Equally important are the cones in \cref{D2.1}.
Note that  although the drift of the diffusion is
piecewise-linear when the control is constant, it becomes
quite complicated under a (stationary) Markov control.
Careful analysis of the drift of the diffusion in \cref{E-drift} on these cones
enables us to obtain the drift inequalities and Foster--Lyapunov equations in
\cref{L2.1,T2.1,T2.2}.
The more specialized results in \cref{S2.6} involve Lyapunov functions
which are sums of two exponentials.

The relation between the prelimit dynamics and the limiting diffusion
can be described as follows.
For a model with Poisson arrivals,
the process $\process{\Hat{X}^n}$  describing the (diffusion-scaled)
total number of jobs in the system
is a controlled Markov process with generator (see \cref{E-hatAg,E-sA})
\begin{equation*}
\widehat\sA^n_z f(\Hat{x}) \,\df\, 
\sum_{i\in\cI} \lambda^n_i \bigl(f(\Hat{x}+n^{-\nicefrac{1}{2}} e_i) - f(\Hat{x})\bigr)
+ \sum_{i\in\cI}\bigl(\mu^n_i z_i +\gamma^n_i q_i(\Hat{x},z)\bigr)
\bigl(f(\Hat{x}-n^{-\nicefrac{1}{2}} e_i) - f(\Hat{x})\bigr)\,.
\end{equation*}
Here, the vector $z=(z_1,\dotsc,z_m)\in \ZZ_+^m$ is the control parameter, with
$z_i$ denoting the number of jobs of class $i$ in service, $\{\lambda^n_i\}_{i\in\cI}$,
$\{\mu^n_i\}_{i\in\cI}$, and $\{\gamma^n_i\}_{i\in\cI}$ are the arrival, service rates,
and abandonment rates,  respectively, and $q=(q_1,\dotsc,q_m)$ is the vector of
queue sizes.
Using the diffusion-scaled variables $\Hat{z}^n$ and $\Hat{q}^n$ defined
in \cref{E-hatx} as
\begin{equation*}
\Hat{z}^n_i \,\df\, \frac{1}{\sqrt{n}}\biggl(z_i - \frac{\lambda^n_i}{\mu^n_i}\biggr)
- \frac{\varrho^n}{m}\,,
\quad\text{and\ \ } \Hat{q}^n_i \,\df\, \frac{q_i(x,z)}{\sqrt{n}}\,,
\end{equation*}
we obtain
\begin{equation}\label{E-sA3}
\begin{aligned}
\widehat\sA^n_z f(\Hat{x}) \,=\, 
&\sum_{i\in\cI} \frac{\lambda_i^n}{n}
\frac{f(\Hat{x}+ n^{-\nicefrac{1}{2}} e_i) -2 f(\Hat{x})
+f(\Hat{x}- n^{-\nicefrac{1}{2}} e_i)}{n^{-1}}\\
&\mspace{200mu}
 + \sum_{i\in\cI}\bigl(\mu_i^n\tfrac{\varrho^n}{m}
+\mu^n_i \Hat{z}_i^n +\gamma^n_i \Hat{q}_i^n\bigr)
\frac{f(\Hat{x}- n^{-\nicefrac{1}{2}} e_i) - f(\Hat{x})}{n^{-\nicefrac{1}{2}}}\,.
\end{aligned}
\end{equation}
As shown in \cref{PT3.1K}, for any work-conserving job allocation $z\in \ZZ_+^m$,
 there exists a vector $u\in\varDelta$ such that
$\Hat{z}^n = \Hat{x} - \langle e, \Hat{x}\rangle^+ u$, and
$\Hat{q}^n=\langle e, \Hat{x}\rangle^+ u$.
Using these identities in \cref{E-sA3}, and letting $n\to\infty$, we obtain the
generator of the controlled diffusion in \cref{E-sde}
(see also \cref{R2.3}).
There is some difficulty though with translating the Foster--Lyapunov equation
for the diffusion into an analogous equation for the operator $\widehat\sA^n_z$.
This is because, whereas $\Hat{z}^n$ is of order $\sqrt n$, the queue sizes $\Hat{q}^n$
are not bounded.
We circumvent this problem by establishing
drift inequalities in \cref{L2.1} for the truncated drift given in \cref{Ebc}.
This facilitates using the same Lyapunov function
for the stability analysis of the diffusion and the prelimit, and consequently, 
we have a unified approach to the problem. 

When studying the diffusion-scaled model with renewal arrivals,
the Lyapunov function has to be augmented to account for the age processes
(see \cref{E3.2-Lyap}).
The analysis is more intricate in this case, and deriving the Foster--Lyapunov
equations in \cref{T3.1} requires extra care.

\begin{remark}\label{R2.3}
Note that if we let $\zeta =\frac{\varrho}{m} e + M^{-1}\ell$,
with $\varrho$ as in \cref{E-varrho},
then a mere translation of the origin of the form $\Tilde X_t = X_t - \zeta$
results in a diffusion of the form \cref{E-sde} with the constant
$\ell$ in the drift taking the form $\ell=- \frac{\varrho}{m} Me$. 
Therefore, without loss of generality, we assume throughout the paper that
the drift in \cref{E-drift} takes the form
\begin{equation}\label{E-drift2}
b(x,u) \,=\, -\frac{\varrho}{m}Me
-M(x-\langle  e,x\rangle^+u)-\langle  e,x\rangle^+\varGamma u\,.
\end{equation}
\end{remark}

For $f\in \Cc^2(\Rd)$ and $u\in\varDelta$, we define
$a(x)=\bigl(a^{ij}(x)\bigr)_{1\le i,j\le m}\df \upsigma(x)\upsigma(x)\transp$, and
\begin{equation}\label{E-Lg}
\Lg_uf(x) \,=\,
\frac{1}{2}\trace\bigl(a(x)\nabla^2f(x)\bigr)
+\bigl\langle b(x,u),\nabla f(x)\bigr\rangle\,,
\end{equation}
with $\nabla^2f$ denoting the Hessian of $f$.

\subsection{Preliminaries.}\label{S2.4}
We start with two definitions.

\begin{definition}\label{D2.1}
For $\delta\in[0,1]$, we define the following cones 
\begin{align*}
\cK_{\delta}^+ &\,\df\, \bigl\{x\in\Rd\,\colon \langle  e,x\rangle
\ge \delta \norm{x}^{}_1\bigr\}\,,\\[5pt]
\cK_{\delta}^- &\,\df\, \bigl\{x\in\Rd\,\colon \langle  e,x\rangle
\le - \delta \norm{x}^{}_1\bigr\}.
\end{align*}
Note that $\cK_0^+$ ($\cK_0^-$) is the nonnegative (nonpositive) canonical half-space,
and $\cK_1^+$ ($\cK_1^-$) is the nonnegative (nonpositive) closed orthant.
Also we have the following identities:
\begin{equation}\label{E-idcone}
\langle  e,x^+\rangle \,=\, \frac{1\pm\delta}{2}\,\norm{x}^{}_1\,,\qquad
\langle  e,x^-\rangle \,=\, \frac{1\mp\delta}{2}\,\norm{x}^{}_1
\qquad\text{for\ \ } x\in\partial\cK_\delta^\pm\,,\ \delta\in[0,1]\,.
\end{equation} 

We fix some convex function $\psi\in\Cc^2(\RR)$ with the property that
$\psi(t)$ is constant for $t\le-1$, and $\psi(t)=t$ for $t\ge0$.
This is defined by 
\begin{equation*}
\psi(t) \df\,
\begin{cases}
-\frac{1}{2}, & t\le-1\,,\\[5pt]
(t+1)^3 -\frac{1}{2}(t+1)^4-\frac{1}{2}\,, & t \in [-1,0]\,,\\[5pt]
t\,, & t\ge 0\,.
\end{cases}
\end{equation*}
For $\veo>0$ we define
\begin{equation*}
\psi_\veo(t) \,\df\, \psi (\veo t)\,,
\end{equation*}
Thus $\psi_\veo(t)=\veo t$ for $t\ge 0$.
A simple calculation also shows that $\psi''_\veo(t)\le \frac{3}{2}\veo^2$.
\end{definition}

\smallskip
\begin{definition}\label{D2.2}
We let $\beta_i\df\frac{\gamma_i}{\mu_i}$ for
$i\in\cI$.
With $\theta$ and $\veo$ positive constants, we define
\begin{equation}\label{E-Psi}
\begin{aligned}
\Psi_\veo(x) \,\df\, \sum_{i\in\cI} \frac{\psi_\veo(x_i)}{\mu_i}\,,
&\quad \Psi(x) \,\df\, \sum_{i\in\cI} \frac{\psi(x_i)}{\mu_i}\,, \\[5pt]
\text{and\ \ } \Psi^*_{\veo,\theta}(x) \,\df\,& \veo\theta \Psi(-x)+ \Psi_\veo(x)\,.
\end{aligned}
\end{equation}
\end{definition}

The function $\Psi$ plays a fundamental role in our analysis.
The quantity $\Psi(x^+)$ represents of course the total workload,
while $\Psi(x^-)$ is a measure of idleness.
These functions, without the smooth cutoff part, are also used
in \cite{GS-12} as test functions to estimate the tails of
the invariant distribution of the prelimit diffusion-scaled processes.

The function $\Psi^*_{\veo,\theta}$
``follows'' the norm $\norm{\,\cdot\,}^{}_1$, in the sense that
\begin{equation}\label{E-Id01}
\veo \tfrac{1\wedge\theta}{\mu\mx}\, \norm{x}^{}_1 - \tfrac{m}{2}
\,\le\, \Psi^*_{\veo,\theta}(x)\,\le\,\veo \tfrac{1\vee\theta}{\mu\mn}\, \norm{x}^{}_1
\,.
\end{equation}
We also have $\psi'(-\nicefrac{1}{2}) = \nicefrac{1}{2}$,
from which we obtain
\begin{equation}\label{E-Id02}
\sum_{i\in\cI} \psi_\veo'(x_i) x_i \,\ge\, \veo\norm{x^+}^{}_1 - \frac{m}{2}\,,
\quad\text{and\ \ }
-\sum_{i\in\cI} \psi'(-x_i) x_i \,\ge\, \norm{x^-}^{}_1 - \frac{m}{2}\,.
\end{equation}
Note also that
\begin{equation}\label{E-Id03}
-\veo\sum_{i\in\cI} \psi'(-x_i) x_i \,\le\, \veo \langle e,x\rangle
\,\le\, \sum_{i\in\cI} \psi_\veo'(x_i) x_i\,.
\end{equation}

Using the parameter $\beta_i$ in \cref{D2.2} and \cref{E-drift2}, we write
the following identities, which we use frequently in the rest of the paper.
\begin{subequations}
\begin{align}
&\begin{aligned}
\bigl\langle \nabla \Psi_\veo(x), b(x,u)\bigr\rangle &\,=\,
-\frac{\varrho}{m} \sum_{i\in\cI} \psi_\veo'(x_i)
- \sum_{i\in\cI} \psi_\veo'(x_i)x_i
+ \langle e,x\rangle^+ \sum_{i\in\cI} \psi_\veo'(x_i)(1-\beta_i)^+ u_i \\
&\mspace{200mu}- \langle e,x\rangle^+ \sum_{i\in\cI} \psi_\veo'(x_i)(\beta_i-1)^+ u_i \,,
\label{E-Id04a}
\end{aligned}\\[5pt]
&\begin{aligned}
\bigl\langle \nabla \Psi(-x), b(x,u)\bigr\rangle &\,=\,
\frac{\varrho}{m} \sum_{i\in\cI} \psi'(-x_i) + \sum_{i\in\cI} \psi'(-x_i)x_i
- \langle e,x\rangle^+ \sum_{i\in\cI} \psi'(-x_i)(1-\beta_i)^+ u_i \\
&\mspace{200mu}+
 \langle e,x\rangle^+ \sum_{i\in\cI} \psi'(-x_i)(\beta_i-1)^+ u_i \,.\label{E-Id04b}
\end{aligned}
\end{align}
\end{subequations}

\subsection{Main results on uniform exponential ergodicity.}\label{S2.5}

The following lemma presents some important drift inequalities
which are used frequently throughout the paper.
Recall the definitions in \cref{E-Psi}.
In order to apply the drift inequalities for the diffusion
to the prelimit in \cref{S3}, we often need to truncate
the diffusion-scaled queueing processes.
To prepare for this, we present a more general version
of these inequalities than what is needed in this section.

For a constant $c\in[1,\infty]$ we define
$b_c(x,u) \df \bigl(b_c^1(x,u),\dotsc,b_c^m(x,u)\bigr)\transp$, with
\begin{equation}\label{Ebc}
b_c^i(x,u) \,\df\, -\frac{\varrho\mu_i}{m}
- \mu_i \bigl(x_i-\langle e,x\rangle^+ u_i\bigr) 
-\gamma_i\langle e,x\rangle^+ u_i \,\Ind_{\{x_i\le c\}}\,.
\end{equation}

\begin{lemma}\label{L2.1}
Assume $\varrho>0$, and let $\theta$ be a positive constant
satisfying
\begin{equation}\label{EL2.1A}
\theta(\beta\mx-1)^+ \,\le\, 1\,.
\end{equation}
Then, the function
\begin{equation}\label{EL2.1B}
V(x)\,\df\, \exp\big(\Psi^*_{\epsilon,\theta}(x)\big) 
\,=\, \exp\bigl(\veo\theta\Psi(-x) + \Psi_\veo(x)\bigr)
\end{equation}
satisfies, for any constant $c\in[1,\infty]$, 
\begin{equation}\label{EL2.1C}
\bigl\langle \nabla V (x), b_c(x,u)\bigr\rangle \,\le\,
\begin{cases}
\veo\Bigl(\theta\varrho +
\frac{m}{2\veo}(1+\veo\theta) - (\theta\wedge1) \norm{x}^{}_1\Bigr) V(x)
&\forall\,x\in\cK_0^-\,,\\[5pt]
- \veo\bigl(\frac{\varrho}{m}-\theta\varrho-\theta\frac{m}{2}
+\theta\norm{x^-}{}_1\bigr)\,V(x) &\forall\,(x,u)\in\cK_0^+\times\varDelta\,.
\end{cases}
\end{equation}
\end{lemma}

\proof{Proof.}
The bound on $\cK_0^-$ follows by first multiplying \cref{E-Id04b} by
$\veo\theta$, then adding this equation to \cref{E-Id04a}, and using \cref{E-Id02}.

We proceed to derive the stated bound on $\cK_0^+\times\varDelta$.
Note that
\begin{equation*}
\bigl\langle \nabla \Psi(-x), b_c(x,u)\bigr\rangle
\,=\,\bigl\langle \nabla \Psi(-x), b(x,u)\bigr\rangle\,,
\end{equation*}
that is, it is equal to the right-hand side of \cref{E-Id04b} for any $c\ge1$,
since $\psi'(-r)=0$ for $r\ge1$.
We write
\begin{equation}\label{PL2.1A}
\begin{aligned}
\bigl\langle \nabla \Psi_\veo(x), b_c(x,u)\bigr\rangle &\,=\,
-\frac{\varrho}{m} \sum_{i\in\cI} \psi_\veo'(x_i)
- \sum_{i\in\cI} \psi_\veo'(x_i)\bigl(x_i-\langle e,x\rangle^+ u_i\bigr)\\
&\mspace{250mu}- \langle e,x\rangle^+ \sum_{i\in\cI} \psi_\veo'(x_i)\beta_i u_i
\,\Ind_{\{x_i\le c\}} \,.
\end{aligned}
\end{equation}
It holds that
\begin{equation}\label{PL2.1B}
-\frac{\varrho}{m} \sum_{i\in\cI} \psi_\veo'(x_i)
\,\le\, -\veo\frac{\varrho}{m}\quad\text{on\ }\cK_0^+\,.
\end{equation}
Also, by \cref{E-Id02}, we have
\begin{equation}\label{PL2.1C}
\begin{aligned}
\theta\veo \frac{\varrho}{m} \sum_{i\in\cI} \psi'(-x_i)
+ \theta\veo \sum_{i\in\cI} \psi'(-x_i)x_i
&\,\le\, \veo\theta\varrho + \veo\theta \frac{m}{2} - \veo\theta\norm{x^-}{}_1
\quad\text{on\ }\Rd\,,
\end{aligned}
\end{equation}
and
\begin{equation}\label{PL2.1D}
\sum_{i\in\cI} \psi_\veo'(x_i)\bigl(x_i-\langle e,x\rangle^+ u_i\bigr)
\,\ge\,0 \quad\text{for\ }x\in\cK_0^+
\end{equation}
by \cref{E-Id03}.
Thus, if $\beta\mx\le1$, then it is
clear from \cref{E-Id04a} and \cref{PL2.1A,PL2.1B,PL2.1C,PL2.1D}, that
\cref{EL2.1C} holds for any positive $\veo$ and $\theta$.

Next suppose $\beta\mx>1$.
We proceed by carefully comparing the terms in \cref{E-Id04b,PL2.1A}.
Define
\begin{equation}\label{E-hcI}
\widehat\cI\,\df\,\{i\in\cI\colon \gamma_i>\mu_i\}\,,\quad
\widehat\cI_+(x)\,\df\,\{i\in\widehat\cI\,\colon x_i\ge0\}\,,\quad\text{and\ \ }
\widehat\cI_-(x)\,\df\,\{i\in\widehat\cI\,\colon x_i<0\}\,.
\end{equation}
Since $\theta(\beta\mx-1)^+\le 1$, and
$\beta_i>1$ on $\widehat\cI$, we have
\begin{equation*}
\veo\theta\sum_{i\,\in\,\widehat\cI_+(x)} \psi'(-x_i)(\beta_i-1)^+ u_i
\,\le\,\sum_{i\,\in\,\widehat\cI_+(x)} \psi_\veo'(x_i)\beta_i u_i\,\Ind_{\{x_i\le c\}}\,,
\end{equation*}
which we combine with
\begin{equation*}
\sum_{i\,\in\,\widehat\cI_-(x)}\psi_\veo'(x_i) u_i
-\sum_{i\,\in\,\widehat\cI_-(x)}\psi_\veo'(x_i)\beta_i u_i\,\Ind_{\{x_i\le c\}}
\,\le\,0\,,
\end{equation*}
to write
\begin{equation}\label{PL2.1E}
\veo\theta\sum_{i\,\in\,\widehat\cI_+(x)} \psi'(-x_i)(\beta_i-1)^+ u_i
+ \sum_{i\,\in\,\widehat\cI_-(x)}\psi_\veo'(x_i) u_i
- \sum_{i\,\in\,\widehat\cI}\psi_\veo'(x_i)\beta_i u_i\,\Ind_{\{x_i\le c\}}
\,\le\, 0\,.
\end{equation}
By the definitions in \cref{E-hcI}, we have the identity
\begin{equation}\label{PL2.1F}
\veo\sum_{i\,\in\,\widehat\cI_-(x)} \psi'(-x_i) u_i\,=\,
\veo\sum_{i\in\widehat\cI} u_i
- \sum_{i\,\in\,\widehat\cI_+(x)} \psi_\veo'(x_i) u_i\,.
\end{equation}
Using again the fact that
$\psi'(-r)=0$ for $r\ge1$, we obtain
\begin{equation}\label{PL2.1G}
\begin{aligned}
\veo\theta\langle e,x\rangle^+ \sum_{i\,\in\,\widehat\cI_-(x)} 
\psi'(-x_i) (\beta_i-1)^+ u_i
&\,\le\, \veo\theta (\beta\mx-1)^+ \langle e,x\rangle^+ 
 \sum_{i\,\in\,\widehat\cI_-(x)} \psi'(-x_i) u_i \\[3pt]
&\,\le\, \veo\langle e,x\rangle^+ \sum_{i\,\in\,\widehat\cI_-(x)} \psi'(-x_i) u_i\\[3pt]
&\,\le\, \veo\langle e,x\rangle^+ \sum_{i\in\widehat\cI} u_i
- \langle e,x\rangle^+ \sum_{i\,\in\,\widehat\cI_+(x)} \psi_\veo'(x_i) u_i
\end{aligned}
\end{equation}
for all $(x,u)\in\cK_0^+\times\varDelta$.
In the second inequality of \cref{PL2.1G} we used the fact
that $\theta(\beta\mx-1)^+\le 1$, and in the third we used \cref{PL2.1F}.
Multiplying \cref{PL2.1E} by by $\langle e,x\rangle^+$, and adding it to
\cref{PL2.1G}, and then combining the resulting sum with the inequality
\begin{equation*}
\langle e,x\rangle^+\sum_{i\,\in\,\widehat\cI\comp}\psi_\veo'(x_i) u_i
-\veo \langle e,x\rangle^+ \sum_{i\,\in\,\widehat\cI\comp} u_i \,\le\, 0\,,
\end{equation*}
where $\widehat\cI\comp$ denote the complement of $\widehat\cI$ with respect to $\cI$,
we obtain 
\begin{equation}\label{PL2.1J}
\begin{aligned}
&\veo\theta\langle e,x\rangle^+ \sum_{i\,\in\,\widehat\cI} \psi'(-x_i) (\beta_i-1)^+ u_i
- \veo \langle e,x\rangle^+
+ \sum_{i\in\cI} \psi_\veo'(x_i) \langle e,x\rangle^+ u_i\\
&\mspace{250mu}
- \langle e,x\rangle^+
\sum_{i\,\in\,\widehat\cI}\beta_i\psi_\veo'(x_i) u_i\,\Ind_{\{x_i\le c\}}
 \,\le\,0\,.
\end{aligned}
\end{equation}
Replacing the term $- \veo \langle e,x\rangle^+$ in \cref{PL2.1J}
with $-\sum_{i\in\cI} \psi_\veo'(x_i)x_i$ preserves this inequality
by \cref{E-Id03}.
Thus, by \cref{E-Id04b,PL2.1A,PL2.1B,PL2.1C,PL2.1J}, we obtain\begin{equation*}
\bigl\langle \nabla \Psi^*_{\veo,\theta}(x), b_c(x,u)\bigr\rangle \,\le\,
- \veo\frac{\varrho}{m}+\veo\theta\varrho+\veo\theta\frac{m}{2}
- \veo\theta\norm{x^-}{}_1
\qquad\forall\,(x,u)\in\cK_0^+\times\varDelta\,,
\end{equation*}
from which the second bound in \cref{EL2.1C} follows.
This completes the proof.
\Halmos
\endproof

Recall the definitions in \cref{E-Lg,E-fnorm}.
Also recall that $\uppi_v$ denotes
the invariant probability measure of the process
governed by \cref{E-sde} for a control $v\in\Usm$, under which
$\process{X}$ is 
positive recurrent.

\begin{theorem}\label{T2.1}
Assume that $\varrho>0$, and in addition to \cref{EL2.1A}, let
\begin{equation}\label{ET2.1A}
0\,<\,\theta \,\le\, \frac{\varrho}{3m(2\varrho+m)}\,.
\end{equation}
Then the following hold:
\begin{enumerate}
\item[\ttup{a}]
There exists $\veo_0>0$,
such that for each $\veo\le\veo_0$, the function $V$ in \cref{EL2.1B},
satisfies the Foster--Lyapunov equation 
\begin{equation}\label{ET2.1B}
\Lg_u V(x) \,\le\, \kappa_0 -
\veo\Bigl(\frac{\varrho}{2m}+\theta\norm{x^-}{}_1\Bigr) V(x)
\qquad\forall (x,u)\in\Rd\times\Delta\,,
\end{equation}
for some positive constant
$\kappa_0$ which depends only on $\veo$ and $\theta$.
In particular, the process $\process{X}$ is positive recurrent under
any control $v\in\Usm$, and
\begin{equation}\label{ET2.1C}
\int_\Rd V(x)\,\uppi_v(\D{x}) \,\le\, \frac{2m}{\veo \varrho }\kappa_0\,.
\end{equation}
\item[\ttup{b}]
There exist positive constants
$\gamma$ and $C_\gamma$ such that
\begin{equation}\label{ET2.1D}
\bnorm{P^v_t(x,\cdot\,)-\uppi_v(\cdot)\,}_{V}\,\le\,
C_\gamma V(x)\, \E^{-\gamma t}\qquad\forall (t,x)\in\RR_+\times\Rd\,,\ 
\forall\, v\in\Usm\,.
\end{equation}
\end{enumerate}
\end{theorem}

\proof{Proof.}
Recall the definitions in \cref{E-Psi}, and also define
\begin{equation*}
\psi^*_{\veo,\theta}(t)\,\df\, \veo\theta\psi(-t) + \psi_\veo(t)\,,
\quad t\in\RR\,.
\end{equation*}
Write the diffusion matrix as
$\upsigma=\diag\bigl(2\Tilde\lambda_1,\dotsc,2\Tilde\lambda_m\bigr)^{\nicefrac{1}{2}}$.
For the queueing network,
$\Tilde\lambda_i=\frac{1}{2}\lambda_i(1+c_{a,i}^2)$,
where $\lambda_i$ is the arrival rate (of the fluid limit), and
$c_{a,i}^2$ is squared coefficient of variation of the renewal arrival process
(see \cref{S3.1}).
In the case of a system with Poisson arrivals, $\Tilde\lambda_i=\lambda_i$,
as in \cref{E-sA3}.
See \cref{S3.1} for the definition of these parameters.
We have
\begin{equation*}
\frac{1}{2}\,\trace\bigl(a\nabla^2V(x)\bigr)
\,=\, \Biggl(\sum_{i\in\cI} \frac{\Tilde\lambda_i}{\mu_i}
(\psi^*_{\veo,\theta})''(x_i)
+ \sum_{i\in\cI} \frac{\Tilde\lambda_i}{\mu_i^2}\bigl[(\psi^*_{\veo,\theta})'(x_i)\bigr]^2
\Biggr) V(x)
\qquad\forall\,x\in\Rd\,.
\end{equation*} 
Recall that $\psi_\veo''\le \frac{3}{2}\veo^2$.
Therefore, since also $\psi_\veo'\le \veo$, $\theta\le1$,
and $\sum_i \frac{\Tilde\lambda_i}{\mu_i}=1$ (see \cref{E-HWpara}), we obtain
\begin{equation}\label{PT2.1A}
\frac{1}{2}\,\trace\bigl(a\nabla^2V(x)\bigr)
\,\le\, \veo\bigl(\tfrac{3}{2}(\veo+\theta)+\veo\Bar{C}
\bigr)
V(x)\,,\qquad
\text{with\ \ } \Bar{C}\,\df\,
\sum_{i\in\cI} \frac{\Tilde\lambda_i}{\mu_i^2}\,.
\end{equation}
We also have $\theta\varrho+\theta\frac{m}{2}\le\frac{\varrho}{6m}$,
and $\frac{3}{2}\theta\le\frac{\varrho}{4m}$ by \cref{ET2.1A}.
Thus \cref{ET2.1B} follows from \cref{EL2.1C} by selecting
$\veo<\frac{\varrho}{6m(3+ 2\Bar{C})}$,
while \cref{ET2.1C} follows by \cref{ET2.1B} and It\^o's formula
in the usual manner.

We now turn to part (b).
Write \cref{ET2.1B} as
\begin{equation}\label{PT2.1B}
\Lg_u V(x) \,\le\, \kappa_0 - \kappa_1 V(x)\,,
\end{equation}
We follow the proof of \cite[Theorem~6.1]{MT-III-93} which uses
a $\delta$-skeleton chain $\{X_{\delta n}\}_{n\in\NN}$.
Note that we can use any $\delta>0$, because $P_\delta^v(x,B)>0$ for any
set $B$ with positive
Lebesgue measure. Thus, for simplicity, we use $\delta=1$.
Then, with $t=n+s$, $s\in[0,1)$, we have
\begin{equation}\label{PT2.1C}
\begin{aligned}
\bnorm{P^v_t(x,\cdot) - \uppi_v(\cdot)}_V &\,=\,
\sup_{\substack{g\in\cB(\Rd), \; \abs{g}\le V}}\;\,
\babss{\int_\Rd P^v_{n+s}(x,\D{y}) g(y) - \int_\Rd g(y)\uppi_v(\D{y})}\\[5pt]
&\,\le\, \int_\Rd P_s^v(x,\D{y}) \bnorm{P^v_n(y,\cdot) - \uppi_v(\cdot)}_V\,.
\end{aligned}
\end{equation}
Next, we estimate
$\bnorm{P^v_n(y,\cdot) - \uppi_v(\cdot)}_V$ using
\cite[Theorem~2.3]{MeynT-94b}.
Using It\^o's formula and \cref{PT2.1B}, we obtain (see \cite[Lemma~2.5.5]{book})
\begin{equation}\label{PT2.1D}
\int_\Rd P_t^v(x,\D{y})V(y)\,=\,\Exp^v_x\bigl[V(X_t)]\,\le\,
\frac{\kappa_0}{\kappa_1} + \E^{-\kappa_1t} V(x)
\qquad\forall (t,x)\in\RR_+\times\Rd\,,\ 
\forall\, v\in\Usm\,.
\end{equation}
Therefore, with $\sB$ a ball such that
$\frac{\kappa_0}{\kappa_1}\le
\E^{-\frac{\kappa_1}{2}}\bigl(1-\E^{-\frac{\kappa_1}{2}}\bigr) V(x)$ for $x\in\sB^c$,
we have
\begin{equation*}
\int_\Rd P_1^v(x,\D{y})V(y)\,\le\, \E^{-\frac{\kappa_1}{2}} V(x)
+ \frac{\kappa_0}{\kappa_1} \Ind_\sB(x)\,,
\end{equation*}
which establishes equation (14) in \cite{MeynT-94b}.

The inequality in \cref{ET2.1C}
implies that the collection of
invariant probability measures $\{\uppi_v\colon v\in\Usm\}$ is tight.
By the invariance of $\uppi_v$, tightness,
and the Harnack inequality applied to the densities of $\uppi_v$
(see \cite[Lemma~3.2.4\,(b)]{book}), we have
\begin{equation*}
\int_\Rd \uppi_v(\D{y}) P^v_{\nicefrac{1}{2}}(y,\sB) \,=\, \uppi_v(\sB)
\,\ge\, \beta_0\,>\,0
\end{equation*}
for some constant $\beta_0$ independent of $v\in\Usm$.
Using tightness once more, we can select a ball $B_R\supset\sB$ such that
\begin{equation*}
\int_{B_R} \uppi_v(\D{y}) P^v_{\nicefrac{1}{2}}(y,\sB)\,\ge\, \frac{\beta_0}{2}\,.
\end{equation*}
This implies that
$\sup_{y\in B_R}\, P^v_{\nicefrac{1}{2}}(y,\sB)\,\ge\, \frac{\beta_0}{2}$.
We now employ the parabolic Harnack inequality for
operators in nondivergence form \cite[Theorem~4.1]{Gruber-84}
(for a simpler statement which uses the notation in
this paper see \cite[Theorem~4.7]{RVIM}).
The parabolic Harnack inequality asserts that there exists a positive constant
$C_{\mathsf{H}}$ such that
\begin{equation*}
\sup_{y\in B_R}\, P^v_{\nicefrac{1}{2}}(y,\sB) \,\le\,
C_{\mathsf{H}}\,\inf_{y\in B_R}\, P^v_{1}(y,\sB)\qquad\forall\,v\in\Usm\,.
\end{equation*}
Therefore,
$P_1^v(x,\sB)\ge\frac{1}{2}C_{\mathsf{H}}^{-1}\beta_0$ for all $x\in \sB$ and $v\in\Usm$.
Thus, with $\delta_0\df \frac{1}{4}C_{\mathsf{H}}^{-1}\beta_0$  we can write
\begin{equation*}
\eta\,\df\, \inf_{y\in \sB}\,
P_1^v(x,\sB)-\delta_0\,\ge\,\delta_0
\qquad\forall\,x\in \sB\,,\ \forall\,v\in\Usm\,,
\end{equation*}
which establishes \cite[equation (23)]{MeynT-94b}.

As seen then from equations (19)--(20) and
(24)--(25) in \cite[Theorem~2.3]{MeynT-94b} there exist positive constants
$C_0$ and $\gamma$ depending only on $\kappa_0$, $\kappa_1$, $\eta$, and $\delta_0$,
such that
\begin{equation}\label{PT2.1E}
\bnorm{P^v_n(x,\cdot) - \uppi_v(\cdot)}_V \,\le\, C_0 \E^{-\gamma n} V(x)\,.
\end{equation}
Thus, using \cref{PT2.1E} in \cref{PT2.1C}, and applying \cref{PT2.1D} once more,
we obtain \cref{ET2.1D} for a constant $C_\gamma$ independent of $v\in\Usm$.
This completes the proof.
\Halmos\endproof

\smallskip
Throughout the paper we let $K_r$, or $K(r)$, for $r>0$, denote the closed cube
\begin{equation}\label{E-cube}
K_r \,\df\, \{x\in\Rd \colon \norm{x}^{}_{1} \le r\}\,.
\end{equation}
We also let $\Bar\psi_\veo=\psi_\veo+1$ so that the function is strictly positive,
and define $\Bar\Psi$ and $\Bar\Psi_\veo$ analogously to \cref{E-Psi}.


\begin{remark}\label{R2.4}
Assume that $\varrho>0$, and consider the function
\begin{equation}\label{ER2.4A}
\sV(x) \,\df\, \bigl(\veo\theta\Bar\Psi(-x) + \Bar\Psi_\veo(x)\bigr)^p
\end{equation}
for some $p\ge1$.
Then it follows directly from the proofs of \cref{L2.1,T2.1}
that there exist positive constants $\veo$, $\theta$,
$\Bar\kappa_0$, $\Bar\kappa_1$, and a
cube $K\subset\Rd$, depending only on $p$, such that
\begin{equation*}
\Lg_u \sV(x) \,\le\,
\begin{cases}
\Bar\kappa_0\,\Ind_{K}(x) - \Bar\kappa_1 \sV(x)
&\forall\,x\in\cK_0^-\,,\\[5pt]
- p\veo\frac{\varrho}{2m}\,\bigl(\sV(x)\bigr)^{\frac{p-1}{p}}
&\forall\,(x,u)\in\cK_0^+\times\varDelta\,.
\end{cases}
\end{equation*}
\end{remark}

In \cref{T2.2} which follows we do not assume that $\varrho>0$.

\begin{theorem}\label{T2.2}
Assume that $\varGamma>0$.
With $\Bar{C}$ as defined in \cref{PT2.1A}, let
\begin{equation}\label{ET2.2A}
\theta\,=\, 
\frac{(1-\beta\mn)\vee\frac{1}{2}}{\beta\mx}\,,
\quad\text{and\ \ }
\veo_0\,\df\, \frac{1}{2\sqrt{\Bar{C}}}\,
\bigl[\theta\wedge\beta\mn\bigl(\beta\mn\wedge\tfrac{1}{2}\bigr)\bigr]
\frac{(1\wedge\theta)\mu\mn}{(1\vee\theta)^2\mu\mx}\,.
\end{equation}
Then, for any $\veo\le\veo_0$, the function
\begin{equation*}
\widetilde{V}(x)\,\df\,
 \exp\Bigl(\tfrac{1}{2}\big[\Psi^*_{\epsilon,\theta}(x)\big]^2 \Bigr) 
\,=\, \exp\Bigl(\tfrac{1}{2}
\bigl[\veo\theta\Psi(-x) + \Psi_\veo(x)\bigr]^2\Bigr)
\end{equation*}
satisfies the Foster--Lyapunov equation 
\begin{equation}\label{ET2.2B}
\Lg_u \widetilde{V}(x) \,\le\, \Tilde\kappa_0
- \veo^2\bigl[\theta\wedge\beta\mn\bigl(\beta\mn\wedge\tfrac{1}{2}\bigr)\bigr]
\tfrac{1\wedge\theta}{2\mu\mx} \norm{x}^{2}_1\, \widetilde{V}(x)
\qquad\forall (x,u)\in\Rd\times\Delta\,,
\end{equation}
for a positive constant
$\Tilde\kappa_0$ which depends only on $\veo$
and the system parameters.
In particular, the process $X_t$ governed by \cref{E-sde}
is uniformly exponentially ergodic,
and the associated invariant probability measures have sub-Gaussian tails.
\end{theorem}

\proof{Proof.}
We borrow some calculations from the proof of \cref{L2.1}.
Using \cref{PL2.1F}, and scaling this with
the new definition of $\theta$ in \cref{ET2.2A}, we have
\begin{equation}\label{PT2.2A}
\begin{aligned}
&\bigl((1-\beta\mn)\vee\tfrac{1}{2}\bigr)
\Biggl(-\sum_{i\in\cI} \psi_\veo'(x_i)x_i
+ \langle e,x\rangle \sum_{i\in\cI} \psi_\veo'(x_i)(1-\beta_i)^+ u_i\Biggr)\\[5pt]
&\mspace{350mu}
+ \veo\theta\langle e,x\rangle \sum_{i\in\cI} \psi'(-x_i)(\beta_i-1)^+ u_i \,\le\,0\,.
\end{aligned}
\end{equation}
Here $\varrho$ is not necessarily positive, so
by \cref{E-Id02} we have
\begin{equation}\label{PT2.2B}
-\frac{\varrho}{m} \sum_{i\in\cI} \psi_\veo'(x_i)
+\veo\theta\frac{\varrho}{m} \sum_{i\in\cI} \psi'(-x_i)
+ \veo\theta\sum_{i\in\cI} \psi'(-x_i)x_i
\,\le\, \veo\Bigl(\abs{\varrho} + \theta \abs{\varrho}+\theta\frac{m}{2}
- \theta\norm{x^-}^{}_1\Bigr)
\end{equation}
on $\Rd$.
Note that
\begin{equation*}
\langle e,x\rangle \sum_{i\in\cI} \psi_\veo'(x_i)(1-\beta_i)^+ u_i
\,\le\, \norm{x^+}^{}_1 (1-\beta\mn)
\sum_{i\in\cI} \psi_\veo'(x_i) u_i
\,\le\, \norm{x^+}^{}_1 (1-\beta\mn)\,.
\end{equation*}
Thus, using \cref{E-Id02}, we have
\begin{equation}\label{PT2.2C}
\bigl(\beta\mn\wedge\tfrac{1}{2}\bigr) \Biggl(-\sum_{i\in\cI} \psi_\veo'(x_i)x_i
+ \langle e,x\rangle \sum_{i\in\cI} \psi_\veo'(x_i)(1-\beta_i)^+ u_i\Biggr)
\,\le\, \veo \bigl(\beta\mn\wedge\tfrac{1}{2}\bigr)
\bigl(\tfrac{m}{2\veo} - \beta\mn\norm{x^+}^{}_1\bigr)\,.
\end{equation}
Let $\Bar\theta\df \theta\wedge\beta\mn\bigl(\beta\mn\wedge\tfrac{1}{2}\bigr)$.
Adding \cref{PT2.2A,PT2.2B,PT2.2C}, using \cref{E-Id04a,E-Id04b}, 
and also \cref{E-Id01}, we obtain
\begin{equation}\label{PT2.2D}
\begin{aligned}
\Psi^*_{\veo,\theta}(x) \bigl\langle \nabla \Psi^*_{\veo,\theta}(x), b(x,u)\bigr\rangle
&\,\le\, 
\veo \Bigl(\abs{\varrho} + \theta \abs{\varrho}+\theta\tfrac{m}{2}
+\bigl(\beta\mn\wedge\tfrac{1}{2}\bigr) \tfrac{m}{2\veo}
-\Bar\theta \norm{x}^{}_1\bigr)\Bigr) \Psi^*_{\veo,\theta}(x)
\\[5pt]
&\,\le\,
\veo \Bigl(\abs{\varrho} + \theta \abs{\varrho}+\theta\tfrac{m}{2}
+\bigl(\beta\mn\wedge\tfrac{1}{2}\bigr)\tfrac{m}{2\veo}\Bigr)
\tfrac{1\vee\theta}{\mu\mn}\, \norm{x}^{}_1\\[5pt]
&\mspace{200mu}
-\veo^2 \Bar\theta
\bigl(\tfrac{1\wedge\theta}{\mu\mx}\norm{x}^{}_1 - \tfrac{m}{2\veo}\bigr)
\norm{x}^{}_1
\\[5pt]
&\,\le\,\veo \Hat{c}_0
-\veo^2 \Bar\theta
\tfrac{1\wedge\theta}{\mu\mx} \norm{x}^{2}_1
\qquad\forall\,(x,u)\in\cK_0^+\times\varDelta\,,
\end{aligned}
\end{equation}
where
\begin{equation*}
\Hat{c}_0 \,\df\, \Bigl(\abs{\varrho} + \theta \abs{\varrho}+\theta\tfrac{m}{2}
+\bigl(\beta\mn\wedge\tfrac{1}{2}\bigr)\tfrac{m}{2\veo}\Bigr)
\tfrac{1\vee\theta}{\mu\mn} + \tfrac{m}{2} \bar{\theta}\,.
\end{equation*}
It is straightforward to verify that \cref{PT2.2D} is also valid
on $\cK_0^-\times\varDelta$.
Following the proof of \cref{T2.1}, we have
\begin{equation}\label{PT2.2E}
\begin{aligned}
\trace\bigl(a\nabla^2\widetilde{V}(x)\bigr) &\,\le\,
\Bigl[\tfrac{3}{2}\veo(\veo+\theta)\Psi^*_{\veo,\theta}(x)
+ \veo^2(1\vee\theta)^2\Bar{C}\Bigl(1+
\bigl(\Psi^*_{\veo,\theta}(x)\bigr)^2\Bigr)\Bigr]\widetilde{V}(x)\\[5pt]
&\,\le\, \veo\Bigl[\veo(1\vee\theta)^2\Bar{C} + \tfrac{3}{2} \veo(\veo+\theta)
\tfrac{1\vee\theta}{\mu\mn}\, \norm{x}^{}_1
+\veo^3 \Bar{C} \tfrac{(1\vee\theta)^4}{\mu^2_{\mathsf{min}}}\, \norm{x}^{2}_1\Bigr]
\widetilde{V}(x)\,.
\end{aligned}
\end{equation}
Combining \cref{PT2.2D,PT2.2E}, we obtain
\begin{equation*}
\begin{aligned}
\trace\bigl(a\nabla^2\widetilde{V}(x)\bigr)+
\bigl\langle \nabla \widetilde{V}(x), b(x,u)\bigr\rangle
&\,\le\,
\Bigl[\veo^2(1\vee\theta)^2\Bar{C}
+\veo\Bigl( \tfrac{3}{2}(\veo+\theta)
\tfrac{1\vee\theta}{\mu\mn} + \Hat{c}_0\Bigr)
\norm{x}^{}_1\\[5pt]
&\mspace{60mu}
-\veo^2\Bigl(\bar{\theta}\tfrac{1\wedge\theta}{\mu\mx}
-\veo^2 \Bar{C} \tfrac{(1\vee\theta)^4}{\mu^2_{\mathsf{min}}}\Bigr) \norm{x}^{2}_1
\Bigr]\widetilde{V}(x)\,,
\end{aligned}
\end{equation*}
from which the validity of \cref{ET2.2B} on $\cK_0^+\times\varDelta$ follows
by selecting $\veo$ sufficiently small.
Verifying the validity of \cref{ET2.2B} on $\cK_0^-\times\varDelta$,
is simpler, and is a straightforward application of \cref{E-Id01,E-Id02,PT2.2E}.
This finishes the proof.
\Halmos\endproof

\begin{remark}
The counterpart of \cref{R2.4} applies relative to \cref{T2.2}.
In particular, the function $\sV$ in \cref{ER2.4A}
for $p>0$ is a Lyapunov function.
Indeed, there exist positive constants
$\veo$, $\theta$, $\Check\kappa_0$ and $\Check\kappa_1$, and a
cube $K\subset\Rd$, depending only on $p$, such that
\begin{equation*}
\Lg_u \sV(x) \,\le\,
\Check\kappa_0\,\Ind_{K}(x) - \Check\kappa_1 \sV(x)
\qquad\forall\,(x,u)\in\cK_0^+\times\varDelta\,.
\end{equation*}
\end{remark}

\begin{remark}\label{R2.6}
It is worth noting that if $\varGamma>0$, then by choosing
$\theta>0$ as in \cref{ET2.2A}, the function
\begin{equation*}
\Breve{V}(x)\,\df\, \exp\bigl(\eta\theta\Psi(-x)+\eta\Psi(x)\bigr)
\end{equation*}
satisfies
\begin{equation*}
\Lg_u \Breve{V}(x) \,\le\,
\Breve\kappa_0 - \Breve\kappa_1\norm{x}^{}_1 \Breve{V}(x)
\qquad\forall\,(x,u)\in\cK_0^+\times\varDelta\,,
\end{equation*}
for all $\eta>0$, and for some positive constants $\Breve\kappa_0$ and $\Breve\kappa_1$
depending only on $\eta$.
Indeed, using \cref{PT2.2A,E-Id04a,E-Id04b}, we deduce,
with $\Hat\theta\df1-\bigl((1-\beta\mn)\vee\tfrac{1}{2}\bigr)$, that
\begin{equation*}
\begin{aligned}
\frac{1}{\eta\Breve{V}(x)}\bigl\langle \nabla \Breve{V}(x), b(x,u)\bigr\rangle &\,=\,
\frac{\varrho}{m} \sum_{i\in\cI} \bigl(\theta\psi'(-x_i)-\psi'(x_i)\bigr)
- \sum_{i\in\cI} \bigl((1-\Hat\theta)\psi_\veo'(x_i)x_i-\theta\psi_\veo'(-x_i)\bigr)x_i
\\
&\mspace{150mu}+ (1-\Hat\theta)
\langle e,x\rangle^+ \sum_{i\in\cI} \psi'(x_i)(1-\beta\mn)^+ u_i \\[3pt]
&\,\le\,
\frac{m}{2}(1-\Hat\theta+\theta)
+\frac{\varrho}{m} \sum_{i\in\cI} \bigl(\theta\psi'(-x_i)-\psi'(x_i)\bigr)
- \bigl(\beta\mn(1-\Hat\theta)\wedge\theta\bigr)\norm{x}^{}_1\,,
\end{aligned}
\end{equation*}
where we also used \cref{E-Id02,E-Id03}.
The rest is routine.
\end{remark}

\subsection{Results concerning the tail of the invariant distribution.}\label{S2.6}

Gamarnik and Stolyar in \cite{GS-12} conjecture that,
provided $\varrho>0$,
$\exp\bigl(\theta\sum_i x_i^-\bigr)$ is integrable
under an invariant probability measure for all $\theta>0$.
They prove this when $\gamma_i\le\mu_i$ for all $i\in\cI$
\cite[Theorem~2\,(i)]{GS-12}.
The proof is for the diffusion-scaled queueing processes,
and relies on a simple comparison
to a system with infinitely many servers.
For this proof to go through though, it seems necessary that
all $i$ satisfy $\gamma_i\le\mu_i$.
We improve upon this result, by showing that 
$\E^{\theta x_i^-}$ is integrable
under an invariant probability measure for all $\theta>0$,
for any $i$ such that $\gamma_i\le\mu_i$.
Of course this proof applies to the limiting diffusion,
but we show in \cref{S3} how to recover this property for the prelimit in \cref{T3.5}.
The general conjecture remains open.

We need some notation.
We let
\begin{equation}\label{E-cI1}
\cI_1\,\df\,\{i\in\cI\colon \gamma_i\le\mu_i\}\,,
\end{equation}
and for a positive constant $\eta$, 
we define
\begin{equation}\label{E-Phi}
\Phi_1(x)\,\df\, \sum_{i\in\cI_1} \frac{\psi(-x_i)}{\mu_i}\,,
\quad\text{and\ \ }
\cV_1(x)\,\df\,\exp\bigl(\eta\Phi_1(x)\bigr)\,.
\end{equation}

\begin{lemma}\label{L2.2}
Assume that $\varrho>0$.
Let $\eta>0$ be arbitrary, and $V(x)=\exp\bigl(\Psi^*_{\veo_0,\theta}(x)\bigr)$,
with $\veo_0$ as in \cref{T2.1}, and the constant $\theta$ chosen to satisfy
\cref{EL2.1A,ET2.1A}.
Then
\begin{equation*}
\Lg_u (\cV_1+V)(x) \,\le\,
\begin{cases}
\kappa_0\,\Ind_{K}(x) - \kappa_1 \norm{x}^{}_1 \bigl(\cV_1(x)+ V(x)\bigr)
&\forall\,x\in\cK_0^-\,,\\[5pt]
\kappa_0\,\Ind_{K}(x)
- \veo_0\frac{\varrho}{8m}\,\bigl(\cV_1(x)+ V(x)\bigr)
&\forall\,(x,u)\in\cK_0^+\times\varDelta
\end{cases}
\end{equation*}
for some positive constants $\kappa_0$ and $\kappa_1$,
and some cube $K\in\Rd$.
\end{lemma}

\proof{Proof.}
Using \cref{E-Id04a,E-Id04b}, we write
\begin{equation}\label{PL2.2A}
\begin{aligned}
\frac{1}{\cV_1(x)}\,
\bigl\langle \nabla \cV_1(x), b(x,u)\bigr\rangle &\,=\,
\frac{1}{2}\eta\abs{\cI_1}
+\eta\frac{\varrho}{m} \sum_{i\in\cI_1}\psi'(-x_i)
- \eta \biggl(\frac{1}{2}\abs{\cI_1}-\sum_{i\in\cI_1} \psi'(-x_i)x_i\biggr)
\\[5pt]
&\mspace{50mu}
- \eta \langle e,x\rangle^+ \sum_{i\in\cI_1} (1-\beta_i)\,\psi'(-x_i)u_i\,.
\end{aligned}
\end{equation}
Let 
\begin{equation*}
H(x)\,\df\, \trace\bigl(a\nabla^2 \Phi(x)\bigr) + \bigl\langle\nabla\Phi(x),
a\bigl(\nabla \Phi(x)+2\nabla\Psi^*_{\veo,\theta}(x)\bigr)\bigr\rangle\,.
\end{equation*}
Recall the definition in \cref{E-cube}.
It is clear from \cref{E-idcone} that we can select
$\delta\in(0,1)$
and $r>0$ such that
\begin{equation}\label{PL2.2B}
\begin{aligned}
H(x) \cV_1(x) + \bigl(\eta\tfrac{m}{2}+\eta\varrho\bigr)\,\cV_1(x)
&\,\le\, \veo_0\,\frac{\varrho}{4m} V(x)\,,\\[5pt]
\mspace{50mu}\text{and\ \ }\cV_1(x) &\,\le\, V(x)\qquad
\quad\forall\, x\in K_{r}\comp \cap\cK_{\delta}^+\,.
\end{aligned}
\end{equation}
Combining \cref{ET2.1B,PL2.2B}, we obtain
\begin{equation}\label{PL2.2D}
\Lg_u (\cV_1+V)(x) \,\le\, \kappa_0 - \frac{\veo_0}{2}
\Bigl(\frac{\varrho}{4m}+\theta\norm{x^-}{}_1\Bigr) \bigl(\cV_1(x)+ V(x)\bigr)
\qquad\forall\, x\in K_{r}\comp \cap\cK_{\delta}^+\,,
\end{equation}
and all $u\in\varDelta$.
By \cref{PL2.2A}, we have
\begin{equation}\label{PL2.2E}
\frac{1}{\cV_1(x)}\,
\bigl\langle \nabla \cV_1(x), b(x,u)\bigr\rangle
\,\le\, \eta\bigl(\tfrac{m}{2}+\varrho\bigr)
-\eta \sum_{i\in\cI_1}x_i^-\,.
\end{equation}
Consider the set
\begin{equation*}
\sK \,\df\, \biggl\{x\in \cK_0^+\setminus\cK_{\delta}^+\,\colon
\tfrac{1}{2}\eta \sum_{i\in\cI_1}x_i^-\,\le\,
\eta\bigl(\tfrac{m}{2}+\varrho\bigr)
+H(x)+\veo_0\,\frac{\varrho}{4m}\biggr\}\,.
\end{equation*}
Since $H$ is bounded on $\Rd$,
it is clear by the definition of $\sK$ that $\cV_1$ 
and $\Lg_u\cV_1$ are both bounded on $\sK$.
Therefore, since $V$ is coercive on $\sK$, that is, 
$\liminf_{\{\abs{x}\to\infty\,,\, x\in\sK\}} V(x)\to\infty$,
there exists $r_\circ>0$ such that
\begin{equation}\label{PL2.2G}
\babs{\Lg_u \cV_1(x)} \,\le\, 
\veo_0\,\frac{\varrho}{4m} V(x)\,,\quad\text{and\ \ } \cV_1(x)\le\cV(x)
\qquad\forall (x,u)\in(\sK\cap K_{r_\circ}\comp)\times\Delta\,.
\end{equation}
On the other hand, we have
\begin{equation}\label{PL2.2H}
\Lg_u \cV_1(x)\,\le\, -\Bigl(\veo_0\,\frac{\varrho}{4m}+\frac{\eta}{2}
\norm{x^-}^{}_1\Bigr)\cV_1(x)
\quad\forall\, (x,u)\in (\cK_0^+\setminus\cK_{\delta}^+)\cap\sK\comp
\end{equation}
by \cref{PL2.2E}.
\Cref{PL2.2G,PL2.2H}, together with \cref{ET2.1B,PL2.2D},
imply that
\begin{equation*}
\Lg_u (\cV_1+V)(x)
\,\le\, \kappa_0 - \frac{\veo_0}{2}
\Bigl(\frac{\varrho}{4m}+\bigl(\theta\wedge\tfrac{\eta}{2}\bigr)
\norm{x^-}{}_1\Bigr) \bigl(\cV_1(x)+ V(x)\bigr)
\qquad\forall\, x\in K_{r\vee r_\circ}\comp \cap\cK_0^+\,.
\end{equation*}

The estimate on $\cK_0^-$ is straightforward.
Indeed, \cref{PL2.2A} shows that $\cV_1$ satisfies
this estimate, and \cref{ET2.1B} asserts the same for $V$.
This completes the proof.
\Halmos\endproof

The following is immediate from \cref{L2.2}.

\begin{corollary}\label{C2.1}
Suppose $\varrho>0$. Then the function
$\exp\Bigl(\eta \sum_{i\in\cI_1}
\frac{\psi(-x_i)}{\mu_i}\Bigr)$ is integrable under
the invariant distribution for any $\eta>0$.
\end{corollary}

In \cite[Theorem~4\,(i)]{GS-12} it is shown that if $\nu$ is any
limit of the invariant distributions of the diffusion-scaled queueing processes, 
then there exists some $\theta$ such that
$f(x)=\exp\bigl(\theta\sum_i (x_i^-)^2\bigr)$ is integrable under $\nu$.
As is pointed out in \cite{GS-12}, this property holds only at the limit.
The function $f$ is not integrable under the stationary distribution
of the prelimit model.
The proof is rather tedious and is approached via truncations
(see \cite[Proposition~12]{GS-12}).
In what follows, we provide a simple proof of this result, by
showing that this property holds for the limiting diffusion.

Recall the definitions in \cref{E-cI1,E-Phi}.

\begin{theorem}\label{T2.3}
Assume that $\varrho>0$, and let
\begin{equation*}
\Phi_\eta(x)\,\df\, \sum_{i\in\cI_1} \frac{\psi_\eta(- x_i)}{\mu_i}\,,
\qquad \widetilde\cV_\eta(x)\,\df\,
\exp\Bigl(\tfrac{1}{2}\bigl[\eta\Phi_\eta(x)\bigr]^2\Bigr)\,,\quad\text{and\ \ }
V(x) \,\df\, \exp\bigl(\Psi^*_{\veo_0,\theta}(x)\bigr)\,,
\end{equation*}
with $\veo_0$ and $\theta$ chosen as in \cref{L2.2}.
Then there exists $\eta>0$, such that the function
$\cV\df\widetilde\cV_\eta V$ satisfies
\begin{equation*}
\Lg_u \cV(x) \,\le\, c_0 - c_1 \cV(x)\qquad \forall\,(x,u)\in\Rd\times\Act\,.
\end{equation*}
\end{theorem}

\proof{Proof.}
As in \cref{PL2.2A}, we have
\begin{equation}\label{T2.3A}
\begin{aligned}
\bigl\langle \nabla \Phi_\eta(x), b(x,u)\bigr\rangle &\,=\,
\frac{1}{2}\abs{\cI_1}
+\frac{\varrho}{m} \sum_{i\in\cI_1}\psi_\eta'(-x_i)
- \biggl(\frac{1}{2}\abs{\cI_1}-\sum_{i\in\cI_1} \psi_\eta'(-x_i)x_i\biggr)
\\[3pt]
&\mspace{50mu}
- \langle e,x\rangle^+ \sum_{i\in\cI_1} (1-\beta_i)\,\psi_\eta'(-x_i)u_i
\qquad \forall\,(x,u)\in\Rd\times\varDelta\,.
\end{aligned}
\end{equation}
Let
\begin{equation*}
\widetilde{H}_\eta(x)\,\df\, \frac{1}{2}
\trace\Bigl(a\nabla^2 \bigl[\Phi_\eta(x)\bigr]^2\Bigr)
+\frac{1}{2} \Bigl\langle\nabla\bigl[\Phi_\eta(x)\bigr]^2,
a\bigl(\nabla \bigl[\Phi_\eta(x)\bigr]^2+2\nabla\Psi^*_{\veo,\theta}(x)\bigr)
\Bigr\rangle\,.
\end{equation*}
Note that $\widetilde{H}_\eta(x)\le c_0\eta^2 + c_1 \eta^4 \bigl[\Phi_\eta(x)\bigr]^2$
for some positive constants $c_0 + c_1$.
Consider the set
\begin{equation*}
\widetilde\sK \,\df\,\biggl\{x\in \cK_0^+\,\colon
 \eta\sum_{i\in\cI_1}x_i^-\,\le\, 
\tfrac{\abs{\cI_1}}{2}+\eta\varrho
+ \Bigl(\eta^2\widetilde{H}_\eta(x)+\veo_0\,\eta \frac{\varrho}{4m}\Bigr)
 \bigl[\Phi_\eta(x)\bigr]^{-1} \biggr\}\,.
\end{equation*}
It is clear that $\Phi_\eta(x)$ is bounded on this set,
and thus the same applies to $\widetilde{H}_\eta$,
and $\bigl\langle \nabla \Phi_\eta(x), b(x,u)\bigr\rangle$.
Thus we have
\begin{equation}\label{T2.3C}
\begin{aligned}
\sup_{x\in \widetilde\sK}\,\Bigl[\eta^2\widetilde{H}_\eta(x)
+\eta^2\Phi_\eta(x)\bigl\langle \nabla \Phi_\eta(x), b(x,u)\bigr\rangle\Bigr]
\,\xrightarrow[\eta\searrow0]{}\,0\,.
\end{aligned}
\end{equation}
However, \cref{T2.3C,ET2.1B} imply that $\eta$ may be selected small enough
so that
\begin{equation}\label{T2.3D}
\Lg_u \cV(x) \,\le\, \kappa_0 -
\veo\Bigl(\frac{\varrho}{4m}+\theta\norm{x^-}{}_1\Bigr) \cV(x)
\qquad\forall (x,u)\in(\cK_0^+\cap\widetilde\sK)\times\Delta\,.
\end{equation}
On the other hand, by \cref{T2.3A} and the definition of 
$\widetilde\sK$, we have
\begin{equation*}
\eta^2\widetilde{H}_\eta(x)
+\eta^2\Phi_\eta(x)\bigl\langle \nabla \Phi_\eta(x), b(x,u)\bigr\rangle
\,\le\,0\qquad\forall\,(x,u)\in(\cK_0^+\cap\widetilde\sK\comp)\times\Delta\,,
\end{equation*}
which also implies \cref{T2.3D} on $(\cK_0^+\cap\widetilde\sK\comp)\times\Delta$.
Since the bound on $\cK^+$ is clear, this completes the proof.
\Halmos\endproof

\section{Uniform ergodicity of multiclass many-server queues.}\label{S3}

For a detailed description of this model, see \cite{ABP15}.
Here we only review the basic structure which is used for our results.
We consider a sequence of $GI/M/n+M$ queues with $m$ classes of customers,
indexed by $n$, which is the number of servers.
Customers of each class form their own queue and are served in the order
of their arrival.

\subsection{Model and assumptions.}\label{S3.1}

Let $A^n_i$, $i\in\cI=\{1,\dotsc,m\}$, 
denote the arrival process of class-$i$ customers
with arrival rate $\lambda^n_i$.
We assume that $\{A^n_i\}_{i\in\cI}$ are renewal processes defined as follows.
Let $\{R_{ij}\,\colon i\in\cI\,,\,j\in\NN\}$ be a collection
of independent positive random variables such that, for each $i\in\cI$,
$\{R_{ij}\}_{j\in\NN}$ have a common distribution function $F_i$ having
a density $f_i$, mean equal to $1$,
and squared coefficient of variation (SCV) $c_{a,i}^2 \in (0,\infty)$.
Let
\begin{equation*}
h_i(\tau)\,\df\, \frac{f_i(\tau)}{1-F_i(\tau)}\,,\quad\text{and\ \ }
\zeta_i (\tau) \,\df\, \frac{\int_{\tau}^\infty
\bigl(1- F_i(r)\bigr)\,\D{r}}{1- F_i(\tau)}
\end{equation*}
for $\tau\ge0$,
denote the \emph{hazard rate} and the 
\emph{mean residual life} functions for each $i\in \cI$, respectively.
The arrival process $A^n_i$ is then given by
\begin{equation*}
A_i^n(t) \,\df\, \max\,\Bigl\{k\ge 0 \colon \sum_{j=1}^{k} R_{ij}
\le \lambda^n_i t\Bigr\}\,,\quad t\ge0\,,\ i\in\cI\,.
\end{equation*}
 
We assume the following structural hypotheses on the collection
$\{F_i\}_{i\in\cI}$, which are enforced in this subsection without
further mention.

\begin{assumption}\label{A3.1}
The distribution functions $\{F_i\}_{i\in\cI}$ satisfy
$F_i(0)=0$, and have a locally bounded
density $f_i$ with unbounded support.
In addition, the 
mean residual life functions $\{\zeta_i\}_{i\in\cI}$ are bounded.
\end{assumption}
 
The service and patience times are exponentially distributed, with
class-dependent rates, $\mu_i$ and $\gamma_i$, respectively, for class-$i$
customers.
The arrival, service and abandonment processes of each class are mutually independent.

The queueing process (counting the number both in service and in queue for each class)
of the $n^{\text{th}}$ system $X^n= \{X^n(t)\colon t\ge 0\}$
is governed by
\begin{align*}
X^n_i(t)\,=\,X^n_i(0) + A^n_i(t) - 
Y^n_i \left( \mu^n_i\int_0^t Z^n_i(s) \D{s} \right)
- R^n_i \left(\gamma^n_i \int_0^t Q^n_i(s) \D{s} \right)
\end{align*}
for $i \in \cI$ and $t\ge 0$.
Here $Y^{n}_{i}$ and $R^{n}_i$, are
mutually independent rate-$1$ Poisson processes, independent of the initial
conditions $X^{n}_i(0)$ and the arrival processes $A^n_i$, for all $i\in\cI$. 
 Also, $Z^n_i(s)$ and $Q^n_i(s)$ represent the numbers of class-$i$
jobs in service and in queue at time $s$, $s\ge 0$, respectively.

\subsubsection{The Halfin--Whitt regime.}

The parameters satisfy the following limits as $n\to\infty$ for all $i \in \cI$:
\begin{equation}\label{E-HWpara}
\begin{gathered}
\frac{\lambda^n_i}{n} \,\to\,\lambda_i\,>\,0\,,\qquad
\mu^n_i \,\to\,\mu_i\,>\,0\,,\qquad
\gamma^n_i \,\to\,\gamma_i\, \ge\,0\,,\\[5pt]
\frac{\lambda^n_i - n \lambda_i}{\sqrt n} \,\to\,\Hat{\lambda}_i\,,\qquad
{\sqrt n}\,(\mu^n_i - \mu_i) \,\to\,\Hat{\mu}_i\,,\\[5pt]
\rho_i^n\,\df\,\frac{\lambda^n_i}{n\mu^n_i} \,\to\,\rho_i \,\df\,
\frac{\lambda_i}{\mu_i} \,<\, 1\,,
\qquad\sum_{i=1}^{m} \rho_i\,=\,1\,.
\end{gathered}
\end{equation}

The assumptions in \cref{E-HWpara} imply that 
\begin{equation}\label{E-varrhon}
\varrho^n\,\df\,
{\sqrt n} \biggl(1-\sum_{i=1}^{m}\frac{\lambda^n_i}{n\mu^n_i}\biggr)
 \,\to\,\varrho\,\df\,\sum_{i=1}^{m}
\frac{\rho_i \Hat{\mu}_i - \Hat{\lambda}_i}{\mu_i} \,\in \RR\,.
\end{equation}

We define the diffusion-scaled variables by
\begin{equation}\label{E-center}
\begin{gathered}
\Hat{X}^n_i(t) \,=\, \frac{1}{\sqrt{n}}\biggl(X^n_i(t)
- \frac{\lambda^n_i}{\mu^n_i}\biggr)
- \frac{\varrho^n}{m}\,,\quad
\Hat{Z}^n_i(t) \,=\, \frac{1}{\sqrt{n}}\biggl(Z^n_i(t)
- \frac{\lambda^n_i}{\mu^n_i}\biggr)
- \frac{\varrho^n}{m}\,, \\[3pt]
\Hat{Q}^n_i(t)=\frac{1}{\sqrt{n}}Q^n_i(t)\,, \quad \text{and} \quad
\Hat{A}^n_i(t)=\frac{1}{\sqrt{n}}\bigl(A^n_i(t) - \lambda^n_i t \bigr)\,,
\quad i \in \cI\,.
\end{gathered}
\end{equation} 
Then, we obtain the following representation of $\Hat{X}^n_i(t)$:
\begin{equation}\label{E-ds}
\begin{aligned}
\Hat{X}^n_i(t) &\,=\,\Hat{X}^n_i(0) -\frac{\varrho^n\mu^n_i}{m} t 
- \mu^n_i \int_{0}^{t} \Hat{Z}^n_i(s)\,\D{s}
- \gamma^n_i \int_{0}^{t} \Hat{Q}^n_i(s)\,\D{s}\\[5pt]
&\mspace{300mu} + \Hat{A}^n_i(t) - \Hat{M}_{S,i}^{n}(t)
- \Hat{M}_{R,i}^{n}(t)\,, \quad t \ge 0\,,
\end{aligned}
\end{equation}
where 
\begin{align*}
\Hat{M}_{Y,i}^{n}(t) & \,\df\, \frac{1}{\sqrt{n}} \biggl(Y^n_i
\biggl( \mu^n_i\int_0^t Z^n_i(s) \D{s} \biggr)
- \mu^n_i \int_{0}^{t} \Hat{Z}^n_i(s)\,\D{s} \biggr)\,, \\[3pt]
\Hat{M}_{R,i}^{n}(t) & \,\df\, \frac{1}{\sqrt{n}}
\biggl( R^n_i \biggl(\gamma^n_i \int_0^t Q^n_i(s) \D{s} \biggr)
- \gamma^n_i \int_0^t Q^n_i(s) \D{s} \biggr)\,,
\end{align*}
and the last two terms $\Hat{M}_{Y,i}^{n}(t)$ and $\Hat{M}_{R,i}^{n}(t)$ are
square integrable martingales associated with the service and abandonment processes,
respectively. 
The martingales are compensated rate-1 Poisson processes with random time changes,
with respect to the natural filtration, see \cite{ABP15}. 

Note that the diffusion-scaled arrival processes satisfy
\begin{equation*}
\Hat{A}^n \,\Rightarrow\, \diag\bigl(\lambda_1c^2_{a,1}, \dotsc, \lambda_m c^2_{a,m}
\bigr)^{\nicefrac{1}{2}} W \quad\text{in\ }(\mathbb{D}_m, J_1)\quad
\text{as\ }n\to \infty\,,
\end{equation*}
where $W$ is a standard $m$-dimensional Wiener process and $(\mathbb{D}_m, J_1)$
represents the space of c{\`a}dl{\`a}g functions in $\RR^m$ endowed with the
Skorokhod $J_1$ topology. 
Assuming that $\Hat{X}^n(0) \Rightarrow X(0)=x_0$ for a constant $x_0\in \RR^m$, 
it then follows that $\Hat{X}^n \Rightarrow X$ in $(\mathbb{D}_m, J_1)$
as $n\to \infty$, where 
the limit process $X$ satisfies \cref{E-sde} with
$\upsigma(X_t) = \diag\bigl(\lambda_1(1+c^2_{a,1}), \dotsc,
\lambda_m (1+ c^2_{a,m})\bigr)^{\nicefrac{1}{2}}$.
In the case of Poisson arrivals, we have $c_{a,i}^2=1$ and thus
$\upsigma(X_t) = \diag\bigl(2\lambda_1, \dotsc, 2\lambda_m\bigr)^{\nicefrac{1}{2}}$. 

\begin{remark}
Note that this scaling is different from that used in \cite{AMR04, HZ-04, ABP15}, 
where the centering term uses $n\rho_i$ for the processes $X^n_i(t)$ and $Z^n_i(t)$.
Here we use the prelimit parameters $ \nicefrac{\lambda^n_i}{\mu^n_i}$ together with
the ``adjustment" $ \nicefrac{\varrho^n}{m}$, which can be regarded as
the ``reallocation" of the ``safety staffing". 
Recall that when $\varrho^n>0$ (and $\varrho>0$), the condition in \cref{E-varrhon}
is equivalent to the positive square-root safety staffing rule (see \cite{W92}). 
In addition, the diffusion-scaled process $\Hat{X}^n$ converges to the limiting
diffusion $X$ with the drift given in \cref{E-drift2}.
That follows from the standard martingale convergence technique in
\cite{PTW07} using the representation of $\Hat{X}^n$ in \cref{E-ds}. 
\end{remark}

\subsubsection{Scheduling policies.}

We define the space
\begin{equation*}
\cZn(x)\,\df\, \bigl\{z \in \ZZ_+^m\,\colon
z_i\le x_i\,,\ \norm{z}^{}_1 = n\wedge \norm{x}^{}_1\bigr\}\,. 
\end{equation*}
A scheduling policy is called (stationary) Markov if
$Z^{n}(t)= z\bigl(X^n(t),S^n(t)\bigr)$ for some function
$z\colon \ZZ_+^m\times\RR_+^m\to \cZn(x)$,
in which case we identify the policy with the function $z$.
Let $S^n(t)=\bigl(S^n_1(t),\dotsc,S^n_m(t)\bigr)$,
where $S^n_i(t)$ denotes the \emph{age process} for class-$i$ customers,
defined by
\begin{equation*}
S^n_i(t) \,\df\, t - \frac{1}{\lambda^n_i}\sum^{A^n_i(t)}_{j=1} R_{ij}\,,
\qquad t \ge 0\,.
\end{equation*}
Let 
\begin{equation}\label{E-ri}
r^n_i(s_i) \,\df\,
\lambda^n_i \frac{f_i(\lambda^n_i s_i)}{1-F_i(\lambda^n_i s_i)}\,, \quad s_i \ge 0\,,
\end{equation}
denote the \emph{scaled hazard rate} function for the interarrival times of $A^n_i(t)$.

Under a Markov policy, the process $(X^n,S^n)$ is Markov with extended
generator
\begin{equation}\label{E-Ag}
\begin{aligned}
\Ag^n_z g(x,s) &\,\df\, \sum_{i\in\cI} \frac{\partial g(x,s)}{\partial s_i}
+\sum_{i\in\cI} r^n_i(s_i) \bigl(g(x+e_i,s-s_i e_i) - g(x,s)\bigr)\\
&\mspace{200mu}+ \sum_{i\in\cI}\bigl(\mu^n_i z_i +\gamma^n_i q_i(x,z)\bigr)
\bigl(g(x-e_i,s) - g(x,s)\bigr)\,,
\end{aligned}
\end{equation}
for $g\in\Cc_b(\Rd\times\Rd)$ and $(x,s)\in\ZZ_+^m\times\RR_+^m$.
Here, $q_i(x,z) = x_i-z_i$, and
$e_i\in\Rd$ denotes the vector
with the $i^{\text{th}}$ element equal to $1$ and the rest of its elements equal to $0$.

Let
\begin{equation}\label{E-hatx}
\Hat{x}^n_i(x) \,\df\, \frac{1}{\sqrt{n}}\biggl(x_i - \frac{\lambda^n_i}{\mu^n_i}\biggr)
- \frac{\varrho^n}{m}\,,\quad
\Hat{z}^n_i(x) \,\df\, \frac{1}{\sqrt{n}}\biggl(z_i - \frac{\lambda^n_i}{\mu^n_i}\biggr)
- \frac{\varrho^n}{m}\,,
\text{\ and\ \ } \Hat{q}^n_i (x,z) \,\df\, \frac{q_i(x,z)}{\sqrt{n}}\,.
\end{equation}

We let $\sX^n$ denote the state space of the process $\Hat{X}^n$.
This is a countable subset of $\Rd$.
Since $x\mapsto\Hat{x}^n(x)$ is invertible, the set
$\cZn(x)$ can be equivalently written as a function of
$\Hat{x}^n$, and abusing the notation we write this
as $\cZn(\Hat{x}^n)$.
In order to keep the notation simple, we often drop the superscript
$n$ from $\Hat{x}^n$, when this is used to denote a
generic element of $\sX^n$.

\subsection{Results with renewal arrivals.}\label{S3.2}

The first main
result is \cref{T3.1} below
which is the counterpart of \cref{T2.1} for the $n^{\mathrm{th}}$ system.
In order to state this theorem and demonstrate its proof, we need some
additional notation which we introduce next.
 
Let $V$ be the function in \cref{EL2.1B} with
$\mu^n$ replacing $\mu$
in its definition, and
parameters $\veo>0$ and $\theta\in(0,1)$.
Let 
\begin{equation*}
\zeta_i^n(\tau)\,\df\, \zeta_i(\lambda_i^n\tau)\,,\quad \uptau\ge0\,,\ i\in\cI\,.
\end{equation*}
In \cref{T3.1} below, we use the Lyapunov function $\sV^n$ defined by
\begin{equation}\label{E3.2-Lyap}
\sV^n(\Hat{x},s)\,\df\, G^n (\Hat{x},s) + V(\Hat{x})\,, \quad
(\Hat{x},s)\in\sX^n\times\RR_+^m\,,
\end{equation}
with
\begin{equation*}
G^n (\Hat{x},s) \,\df\,
\sum_{i\in\cI} \bigl(1-\zeta_i^n(s_i)\bigr)
\bigl(V(\Hat{x}+ n^{-\nicefrac{1}{2}} e_i)-V(\Hat{x})\bigr)\,,\quad
(\Hat{x},s)\in\sX^n\times\RR_+^m\,.
\end{equation*}
Note that, by \cref{A3.1}, for any fixed $\theta$,
we can choose $\Tilde\veo_0=\Tilde\veo_0(\theta)>0$ small enough so that
\begin{equation}\label{ES3.2A}
\veo \Babss{\sum_{i\in\cI} \frac{1}{\mu_i^n}\bigl(1-\zeta_i^n(s_i)\bigr)
\bigl(\theta\psi'(-y_i)+\psi'(y_i)\bigr)}
\,\le\, \frac{1}{2} \quad\forall\,\veo\le\Tilde\veo_0(\theta)\,,\ 
\forall\,(y,s)\in\Rd\times\RR_+^m\,, \ \forall n \in \NN\,.
\end{equation}
Then, provided $\veo\le\Tilde\veo_0(\theta)$, we have
\begin{equation}\label{EsVbound}
\frac{1}{2}V(y)\,\le\,\sV^n(y,s)\,\le\, \frac{3}{2} V(y)\,.
\end{equation}
We define
\begin{equation}\label{E-Vhatn}
\widehat{V}^n(x)\,\df\, V\bigl(\Hat{x}^n(x)\bigr)\,,
\quad\text{and\ \ }
\widehat{G}^n(x,s)\,\df\,
\sum_{i\in\cI} \bigl(1-\zeta_i^n(s_i)\bigr)
\Bigl(V\bigl(\Hat{x}^n(x+ e_i)\bigr)-V\bigl(\Hat{x}^n(x)\bigr)\Bigr)
\end{equation}
for $x\in\ZZ_+^m$.
Then
the generator $\widehat\Ag^n_z$ of the diffusion-scaled
state process $(\Hat{X}^n, S^n)$ under a policy $z$ takes the form
\begin{equation} \label{E-hatAg} 
\widehat\Ag^n_z V(\Hat{x},s)\,=\,\Ag^n_z\widehat{V}^n(x,s)\,,\quad\text{and\ \ }
\widehat\Ag^n_z G(\Hat{x},s) \,=\, \Ag^n_z \widehat{G}^n(x,s)\,,
\end{equation}
where $\Ag^n_z$ is as defined in \cref{E-Ag}.

We need to introduce some constants used in the results.
First, for a function $f$ on $\Rd$, if we define
\begin{equation*}
\dd f(x;y) \,\df\, f(x+y)-f(x)\,,
\end{equation*}
it then follows by a repeated use of the mean value theorem that there exists a constant
$\widehat{C}_1$
such that
\begin{equation}\label{E-hC1}
\babs{\dd\widehat{V}^n(x\pm e_j;\pm e_i)-\dd\widehat{V}^n(x;\pm e_i)}
\,\le\, \frac{1}{n} \widehat{C}_1 \veo(\veo+\theta) \widehat{V}^n(x)
\qquad\forall\,i,j\in\cI\,,
\end{equation}
and the same bound holds for
$\babs{\dd\widehat{V}^n(x;e_i)+\dd\widehat{V}^n(x;- e_i)}$.
Also, by \cref{A3.1}, \cref{E-ri}, and the convergence of the parameters in
\cref{E-HWpara}, there exists a constant $\widehat{C}_0^n$ depending
on $n$ (implicitly through $\lambda_i^n$), such that
\begin{equation}\label{E-hC0}
\sup_{n\in\NN}\,\max_{i\in\cI}\,
\biggl(\frac{r_i^n(\tau)}{n}\vee \bigl(1+\zeta_i^n(\tau)\bigr)\biggr)
\,\le\,\widehat{C}_0^n\qquad\forall\,\tau\ge0\,.
\end{equation}
We define
\begin{equation}\label{E-tC}
\widetilde{C}_0^n\,\df\,m^2 \widehat{C}_0^n \widehat{C}_1\,,
\qquad \widetilde{C}_1^n \,\df\,
\widehat{C}_1\Biggl(m^2\widehat{C}_0^n \mu_i^n
+m(m-1) \bigl(\widehat{C}_0^n\bigr)^2
+ \sum_{i\in\cI} \frac{\lambda_i^n}{n}\Biggr)\,,
\end{equation}
and
\begin{equation}\label{PT3.1I}
\theta_0(n) \,\df\,
\frac{1}{1+(\beta^n_\mx-1)^+}\wedge
\frac{1}{2\mu^n_{\mathsf{max}} (\widetilde{C}_0^n+\widehat{C}_1)} \wedge
\frac{\varrho^n}{m} \Bigl(m+2\varrho^n+ 4
\bigl(\widetilde{C}_1^n +  m\widehat{C}_1\,\widehat{C}_2^n+ m\widehat{C}_3^n\bigr)
\Bigr)^{-1}.
\end{equation}

Recall $\Tilde\veo_0(\theta)$ in \cref{ES3.2A}.
We are ready to state the first main result of this section.

\begin{theorem}\label{T3.1}
We enforce \cref{A3.1}, and, in addition, we
assume that
the hazard rate functions $\{h_i\}_{i\in\cI}$ are bounded. Suppose $\varrho^n>0$. 
Then there exists a positive constant
$C_0^n(\veo)$, such that
the function $\sV^n$ in \cref{E3.2-Lyap},
with parameters $\theta=\theta_0(n)$ and
any $\veo<\theta_0(n)\wedge\Tilde{\veo}_0(\theta)$,
satisfies
\begin{equation}\label{ET3.1A}
\widehat\Ag^n_z \sV^n(\Hat{x},s) \,\le\,
C_0^n(\veo) - \veo\frac{\varrho^n}{3m} \sV^n(\Hat{x},s) \qquad\forall\,
(\Hat{x},s)\in \sX^n\times\RR_+^m\,,\ \forall\,z\in\cZn(\Hat{x})\,.
\end{equation}
In particular, under any work-conserving stationary Markov policy,
the process $(\Hat X^n,S^n)$ is positive Harris recurrent,
and $V(\Hat{x})$ is integrable under its invariant probability distribution.
\end{theorem}

\begin{remark}
It is clear from the Foster--Lyapunov equation \cref{ET3.1A},
 that the stability result in \cref{T3.1} holds for
all $n\in\NN$ such that $\varrho^n>0$, and the same applies to \cref{T3.3,C3.1}.
We want to emphasize that this is an
important byproduct of the approach in this paper.
One should compare it to \cite[Theorem~2]{GS-12}, where
stability is only stated as an asymptotic property, or in other words,
that it holds for all large enough $n$.

The convergence
of the parameters in \cref{E-HWpara}, implies
that if the limiting value $\varrho=\lim_{n\to\infty}\varrho^n$ is
positive, then $\theta_0(n)$ and $C_0(\veo)$
can be selected independent of $n$, in a manner that
\cref{ET3.1A} holds for all sufficiently large $n$.
Analogous conclusions can be drawn for
\cref{T3.2,T3.3,C3.1} which appear later in this section.

Note also that the difference in the constant multiplying
the Lyapunov function between \cref{ET2.1B,ET3.1A} is only due
to the bound in \cref{EsVbound}.
\end{remark}

\smallskip
For the proof of the \cref{T3.1} we need the following result.

\begin{lemma}\label{L3.1}
With $\widehat\sV^n(x,s) \,\df\, \widehat{G}^n(x,s) + \widehat{V}^n(x)$,
we have the following inequality
\begin{equation}\label{EL3.1A}
\begin{aligned}
\Ag^n_z\, \widehat\sV^n(x,s) &\,\le\,
 \sum_{i\in\cI} \biggl(\frac{\varrho^n \mu^n_i}{m}
+ \mu^n_i\Hat{z}_i + \gamma^n_i \Hat{q}_i(x,z)\biggr)
\sqrt{n}\,\dd\widehat{V}^n(x;-e_i)\\
&\mspace{150mu}+ \veo(\veo+\theta)
\frac{1}{\sqrt n}\widetilde{C}_0^n
\sum_{i\in\cI} \gamma^n_i\Hat{q}_i(x,z)\widehat{V}^n(x)
+ \veo(\veo+\theta) \widetilde{C}_1^n \widehat{V}^n(x)\,,
\end{aligned}
\end{equation}
with $\widetilde{C}_0^n$ and $\widetilde{C}_1^n$ as defined
in \cref{E-tC}.
\end{lemma}

\proof{Proof.}
Recall the definitions in \cref{E-Vhatn}, and note that
\begin{equation*}
\widehat{G}^n_i (x,s) \,=\, \bigl(1-\zeta_i^n(s_i)\bigr)
\dd\widehat{V}^n(x;e_i)\,,
\end{equation*}
with $\widehat{V}^n$ as defined in \cref{E-Vhatn}.
It follows by direct differentiation that
\begin{equation}\label{E-zetaId}
\frac{\D\zeta_i^n(\tau)}{\D{\tau}} - r_i^n(\tau) \zeta_i^n(\tau)
\,=\, -\lambda_i^n\,,
\quad \tau\ge0\,. 
\end{equation}
Thus, using \cref{E-zetaId,E-hC1,E-hC0}, and noting that $\zeta^n_i(0)=1$,
we obtain 
\begin{equation}\label{PL3.1B}
\begin{aligned}
\Ag^n_z \widehat{G}^n_i(x,s) &\,=\,
-\biggl(\frac{\D\zeta_i^n(s_i)}{\D s_i}
 + r_i^n (s_i)\bigl(1-\zeta_i^n(s_i)\bigr)\biggr)
\dd\widehat{V}^n(x;e_i)\\
&\mspace{10mu}+ r_i^n (s_i)\bigl(1-\zeta_i^n(s_i)\bigr) \sum_{j\ne i\,,\,i\in\cI} 
\bigl(\dd\widehat{V}^n(x+e_j;e_i)-\dd\widehat{V}^n(x;e_i)\bigr)\\[3pt]
&\mspace{20mu}
- \bigl(\mu^n_i z_i +\gamma^n_i q_i(x,z)\bigr) \bigl(1-\zeta_i^n(s_i)\bigr)
\sum_{j\in\cI} \bigl(\dd\widehat{V}^n(x-e_j;e_i)-\dd\widehat{V}^n(x;e_i)\bigr)
\\[3pt]
&\,\le\, \bigl(\lambda_i^n-r_i^n(s_i)\bigr) \dd\widehat{V}^n(x;e_i)
+ (m-1) \bigl(\widehat{C}_0^n\bigr)^2
\widehat{C}_1 \veo(\veo+\theta) \widehat{V}^n(x)\\[5pt]
&\mspace{200mu}
+\frac{m}{n} \widehat{C}_0^n \widehat{C}_1 \veo(\veo+\theta)
\bigl(\mu^n_i z_i +\gamma^n_i q_i(x,z)\bigr) \widehat{V}^n(x)\,.
\end{aligned}
\end{equation}

Also, 
\begin{equation}\label{PL3.1C}
\Ag^n_z \widehat{V}^n(x) \,=\,
\sum_{i\in\cI} r_i^n(s_i) \dd\widehat{V}^n(x;e_i)
+ \sum_{i\in\cI}\bigl(\mu^n_i z_i + \gamma^n_i q_i(x,z)\bigr)
\dd\widehat{V}^n(x;-e_i)\,.
\end{equation}
Applying the identities
\begin{equation}\label{PL3.1D}
z_i \,=\, \sqrt{n}\Hat{z_i} +\frac{\lambda^n_i}{\mu^n_i}
+\sqrt{n}\frac{\varrho^n}{m}\,,\quad\text{and}\quad
q_i(x,z) \,=\, \sqrt{n}\Hat{q}_i(x,z)\,, 
\end{equation}
to \cref{PL3.1C}, we obtain
\begin{equation}\label{PL3.1E}
\begin{aligned}
\Ag^n_z\, \widehat{V}^n(x) &\,=\,
\sum_{i\in\cI} \Big(
r_i^n(s_i)\dd\widehat{V}^n(x;e_i) \bm+\lambda_i^n \dd\widehat{V}^n(x;-e_i) \Big)\\
&\mspace{50mu}
+ \sum_{i\in\cI} \biggl(\frac{\varrho^n \mu^n_i}{m}
+ \mu^n_i\Hat{z}_i + \gamma^n_i \Hat{q}_i(x,z)\biggr)
\sqrt{n}\,\dd\widehat{V}^n(x;-e_i)\,.
\end{aligned}
\end{equation}
Combining \cref{PL3.1B,PL3.1E}, and applying once more the estimate in
\cref{E-hC1} and the inequality $\abs{z_i}\le n$, we deduce that
\begin{equation*}
\begin{aligned}
\Ag^n_z\, \widehat\sV^n(x,s) &\,\le\,
\sum_{i\in\cI} \biggl(\frac{\varrho^n \mu^n_i}{m}
+ \mu^n_i\Hat{z}_i + \gamma^n_i \Hat{q}_i(x,z)\biggr)
\sqrt{n}\,\dd\widehat{V}^n(x;-e_i)\\
&\mspace{50mu}+ \veo(\veo+\theta)
\frac{m^2}{\sqrt n} \widehat{C}_0^n \widehat{C}_1 \gamma^n_i
\Hat{q}_i(x,z)\widehat{V}^n(x)\\
&\mspace{100mu}+ \veo(\veo+\theta) \widehat{C}_1\biggl(
m^2\widehat{C}_0^n \mu_i^n
+m(m-1) \bigl(\widehat{C}_0^n\bigr)^2 
+ \sum_{i\in\cI} \frac{\lambda_i^n}{n}
\biggr)\widehat{V}^n(x)\,.
\end{aligned}
\end{equation*}
This completes the proof.
\Halmos\endproof

\proof{Proof of \cref{T3.1}.}

The proof relies on comparing the right hand side of \cref{EL3.1A}
to the drift inequalities in \cref{L2.1}.
First we fix $n\in\NN$, and as done earlier, we suppress the $n$-dependence of
$\Hat{x}_i^n$, $\Hat{z}_i^n$, and $\Hat{q}_i^n$ in the calculations,
in the interest of simplifying the notation.
It is clear from \cref{E-HWpara,E-hatx} that $\Hat{q}_i\ge0$ if
$x_i\ge n$, or equivalently, if
\begin{equation*}
\Hat{x}_i \,\ge\, \vartheta_n\,\df\,\sqrt{n} (1-\rho_i^n)
- \frac{\sqrt{n}}{m}\,
\biggl(1-\sum_{i=1}^{m}\frac{\lambda^n_i}{n\mu^n_i}\biggr)\,\ge\,0\,.
\end{equation*}
If $\veo\le1$, then
$\psi_\veo(x-y)-\psi_\veo(x)\le -\veo\frac{y}{2}$
and $\psi(-x-y)-\psi_(y) \le y$ for all
$x\ge0$ and $y\in[0,1]$ by \cref{D2.1}.
Thus, if $\theta\in(0,\nicefrac{1}{2}]$ and $\veo\in(0,1]$,
then $V(x-y) \le V(x)$ for all
$x\ge0$ and $y\in[0,1]$.
This of course implies, since $\theta_0(n)<\nicefrac{1}{2}$, that 
$\dd \widehat{V}^n (x; -e_i)\le 0$ if $\Hat{x}_i\ge0$.
Thus, if we write
\begin{equation}\label{PT3.1A}
\begin{aligned}
\biggl(\frac{\varrho^n \mu^n_i}{m}
+ \mu^n_i\Hat{z}_i + \gamma^n_i \Hat{q}_i(x,z)\biggr) &
\sqrt{n}\,\dd\widehat{V}^n(x;-e_i)\\[5pt]
&\,=\,\biggl(\frac{\varrho^n \mu^n_i}{m}
+ \mu^n_i\Hat{z}_i + \gamma^n_i \Hat{q}_i(x,z)\,
\Ind_{\{\Hat{x}_i < \vartheta_n\}}\biggr)
\sqrt{n}\,\dd\widehat{V}^n(x;-e_i)\\[5pt]
&\mspace{100mu}+ \gamma^n_i \Hat{q}_i(x,z)\,
\Ind_{\{\Hat{x}_i \ge \vartheta_n\}}
\sqrt{n}\,\dd\widehat{V}^n(x;-e_i)\,,
\end{aligned}
\end{equation}
then the second term on the right-hand side of \cref{PT3.1A} is negative.
It is also clear that
\begin{equation}\label{PT3.1B}
\babss{\frac{\varrho^n \mu^n_i}{m}
+ \mu^n_i\Hat{z}_i + \gamma^n_i \Hat{q}_i(x,z)\,
\Ind_{\{\Hat{x}_i < \vartheta_n\}}} \,\le\, \widehat{C}_2^n \sqrt{n}
\end{equation}
for some constant $\widehat{C}_2^n$ depending on the parameters.

Using the identity
\begin{equation}\label{PT3.1C}
\widehat{V}^n(x\pm e_i) - \widehat{V}^n(x) \mp \partial_{x_i} \widehat{V}^n(x)
\,=\, \int_0^1 (1-t)\, \partial_{x_ix_i} \widehat{V}^n(x\pm t e_i) \,\D{t}\,,
\end{equation}
we deduce that
\begin{equation}\label{PT3.1D}
\babs{\widehat{V}^n(x\pm e_i) - \widehat{V}^n(x) \mp
\partial_{x_i} \widehat{V}^n(x)}
\,\le\, \frac{1}{n}\veo(\veo+\theta)\,\widehat{C}_1\,\widehat{V}^n(x)\,,
\end{equation}
where, we use a common constant to satisfy \cref{E-hC1,PT3.1D}.
Thus, by \cref{PT3.1A,PT3.1B,PT3.1D}, and
using also the identity
\begin{equation*}
\partial_{x_i} \widehat{V}^n(x)\,=\,\frac{1}{\sqrt n}\partial_{\Hat{x}_i} V(\Hat{x})\,,
\end{equation*}
we obtain
\begin{equation}\label{PT3.1F}
\begin{aligned}
\biggl(\frac{\varrho^n \mu^n_i}{m}
+ \mu^n_i\Hat{z}_i +& \gamma^n_i \Hat{q}_i(x,z)\biggr) 
\sqrt{n}\,\dd\widehat{V}^n(x;-e_i)\\[5pt]
&\,=\,-\biggl(\frac{\varrho^n \mu^n_i}{m}
+ \mu^n_i\Hat{z}_i + \gamma^n_i \Hat{q}_i(x,z)\,
\Ind_{\{\Hat{x}_i < \vartheta_n\}}\biggr)\partial_{\Hat{x}_i} V(\Hat{x})\\[5pt]
&\mspace{100mu}+ \gamma^n_i \Hat{q}_i(x,z)\,
\Ind_{\{\Hat{x}_i \ge \vartheta_n\}}
\sqrt{n}\,\dd\widehat{V}^n(x;-e_i)
+ \veo(\veo+\theta) \widehat{C}_1\,\widehat{C}_2^n\, \widehat{V}^n(x)\,.
\end{aligned}
\end{equation}
Similarly, addressing the second term on the right-hand side
of \cref{EL3.1A}, we write
\begin{equation}\label{PT3.1G}
\frac{1}{\sqrt n}\widetilde{C}_0^n \gamma^n_i\Hat{q}_i(x,z)
\,\le\, \widehat{C}_3^n + 
\frac{1}{\sqrt n}\widetilde{C}_0^n \gamma^n_i\Hat{q}_i(x,z)
\Ind_{\{\Hat{x}_i \ge \vartheta_n\}}
\end{equation}
for some constant $\widehat{C}_3^n$.
Using \cref{E-hatAg,PT3.1F,PT3.1G}, we deduce from \cref{EL3.1A} that
\begin{equation}\label{PT3.1H}
\begin{aligned}
\widehat\Ag^n_z\, \sV^n(\Hat{x},s) &\,\le\,
- \sum_{i\in\cI} \biggl(\frac{\varrho^n \mu^n_i}{m}
+ \mu^n_i\Hat{z}_i + \gamma^n_i \Hat{q}_i(\Hat{x},\Hat{z})\,
\Ind_{\{\Hat{x}_i < \vartheta_n\}}\biggr)\partial_{\Hat{x}_i} V(\Hat{x})\\[5pt]
&\mspace{30mu}+ \veo(\veo+\theta)
\bigl(\widetilde{C}_1^n + m\widehat{C}_1\,\widehat{C}_2^n+ m\widehat{C}_3^n\bigr)
V(\Hat{x})\\[5pt]
&\mspace{50mu}+
\sum_{i\in\cI}\biggl(\sqrt{n}\,\dd\widehat{V}^n(x;-e_i) +
\veo(\veo+\theta)
\frac{1}{\sqrt n}\widetilde{C}_0^n \widehat{V}^n(x)\biggr)
 \gamma^n_i\Hat{q}_i(\Hat{x},\Hat{z})
\Ind_{\{\Hat{x}_i \ge \vartheta_n\}}\,,
\end{aligned}
\end{equation}
where we express $\Hat{q}$ as a function of $\Hat{x}$ and $\Hat{z}$, slightly
abusing the notation.

We now turn to the drift inequalities in \cref{L2.1}.
It follows by \cref{EL2.1C} that there exists a constant
and $C^n_0(\veo)$, such that
\begin{equation}\label{PT3.1J}
\begin{aligned}
&\sum_{i\in\cI} \Bigl(-\frac{\varrho^n\mu^n_i}{m}
-\mu^n_i(\Hat{x}_i-\langle  e,\Hat{x}\rangle^+ u_i)
-\gamma^n_i \langle  e,\Hat{x}\rangle^+ u_i\Ind_{\{\Hat{x}_i < \vartheta_n\}}
\Bigr)\partial_{\Hat{x}_i} V(\Hat{x})\\
&\mspace{150mu} + \veo(\veo+\theta)
\bigl(\widetilde{C}_1^n +  m\widehat{C}_1\,\widehat{C}_2^n+ m\widehat{C}_3^n\bigr)
V(\Hat{x})
\,\le \, C^n_0(\veo) - \veo\frac{\varrho^n}{2m} V(\Hat{x})
\end{aligned}
\end{equation}
for all $(\Hat{x},u)\in\Rd\times\varDelta$,
and for all $\veo\in\bigl(0,\theta_0(n)\bigr)$.

Consider the first sum in \cref{PT3.1J}. 
If $\langle e,x\rangle\le n$, then
$\Hat{z} = \Hat{x}$ by work-conservation.
Note also that by the scaling in \cref{E-center} combined with \cref{E-varrhon},
we have
\begin{equation*}
\langle e, \Hat{x}\rangle \,=\, \frac{1}{\sqrt{n}}
\bigl(\langle e,x\rangle-n\bigr)\,.
\end{equation*}
Thus $\langle e, \Hat{x}\rangle>0$ if and only if $\langle e,x\rangle>n$.
Similarly $\langle e, z\rangle =n$ if and only if $\langle e, \Hat{z}\rangle=0$.
On the other hand,
if $\langle e, x-z\rangle >0$, then we can write
$z = x -\langle e,x-z\rangle u$, for some $u\in\varDelta$.
Thus, $\Hat{z} = \Hat{x} - \langle e, \Hat{x} - \Hat{z}\rangle u
= \Hat{x} - \langle e, \Hat{x} \rangle u$,
since $\langle e, \Hat{z}\rangle=0$.
We have thus established that for all $x\in\RR^m_+$, we have
\begin{equation}\label{PT3.1K}
\Hat{z} \,=\, \Hat{x} - \langle e, \Hat{x}\rangle^+ u\,,\quad
\text{and\ \ } \Hat{q}(\Hat{x},\Hat{z})\,=\,\langle e, \Hat{x}\rangle^+ u
\end{equation}
for some $u\in\varDelta$.
It then follows from \cref{PT3.1K}, that the sum of the first
two terms on the right-hand of
\cref{PT3.1H} has the bound on the right-hand side in \cref{PT3.1J}.

Next, consider the last term in \cref{PT3.1H}.
By \cref{PT3.1C,PT3.1D} we have
\begin{equation*}
\dd \widehat{V}^n(x;-e_i) \,\le\, - \partial_{x_i} \widehat{V}^n(x)
+ \frac{1}{n}\veo(\veo+\theta)\,\widehat{C}_1\,\widehat{V}^n(x)\,,
\end{equation*}
and 
$\partial_{x_i} \widehat{V}^n(x) = \frac{\veo}{\sqrt n\mu_i^n} \widehat{V}^n(x)$
when $\Hat{x}_i \ge \vartheta_n$.
Thus
\begin{equation*}
\biggl(\sqrt{n}\,\dd\widehat{V}^n(x;-e_i) +
\veo(\veo+\theta)
\frac{1}{\sqrt n}\widetilde{C}_0^n \widehat{V}^n(x)\biggr)
\Ind_{\{\Hat{x}_i \ge \vartheta_n\}}
\,\le\, - \veo\biggl(\frac{1}{\mu_i^n}
-(\veo+\theta)
\frac{1}{\sqrt n}(\widetilde{C}_0^n+\widehat{C}_1)\biggr)\widehat{V}^n(x)\,,
\end{equation*}
which is negative for all $\veo<\theta=\theta_0(n)$, by the definition
of $\theta_0$ in \cref{PT3.1I}.
Thus, in view of \cref{EsVbound}, we have established the Foster--Lyapunov equation
\cref{ET3.1A} as claimed.

The remaining conclusions of the theorem are straightforward, in view of
the fact that $\process{S^n}$ is positive Harris recurrent, as shown
in \cite{Takis-99}.
Note that since $\veo<\Tilde\veo_0(\theta)$,
then \cref{EsVbound} implies that $\sV^n$ is bounded from below in $\Rd\times\RR_+^m$.
\Halmos\endproof

In the theorem that follows we assume strictly positive abandonment rates
for all classes, and we use the Lyapunov function
\begin{equation}\label{ET3.2A}
\cV^n(\Hat{x}^n,s)\,\df\,
\sum_{i\in\cI} \bigl(1-\zeta_i^n(s_i)\bigr)
\bigl(\varphi^n(\Hat{x}^n_i+n^{-\nicefrac{1}{2}})-\varphi^n(\Hat{x}^n_i)\bigr)
+ \sum_{i\in\cI} \frac{\varphi^n(\Hat{x}^n_i)}{\mu_i^n}\,,
\end{equation}
with 
\begin{equation*}
\varphi^n(y)\,\df\,\Tilde\veo_0(\theta^{n})\theta^n \psi(-y)
+ \Tilde\veo_0(\theta^{n})\psi(y)\,,
\quad y\in\RR\,,
\end{equation*}
$\Tilde\veo_0$ as in \cref{ES3.2A},
and
\begin{equation}\label{ET3.2B}
\theta^n\,=\, 
1\wedge\frac{(1-\beta^n_{\mathsf{min}})\vee\frac{1}{2}}{\beta^n_{\mathsf{max}}}\,,
\quad\beta^n_i\df \frac{\gamma_i^n}{\mu_i^n}\,,\quad i\in\cI\,.
\end{equation}

\begin{theorem}\label{T3.2}
Grant \cref{A3.1}, and, in addition,
assume that hazard rate functions $\{h_i\}_{i\in\cI}$ are locally bounded.
Suppose $\gamma^n_i>0$ for all $i\in\cI$. 
Then there exist
positive constants $c_0(n)$ and $c_1(n)$, depending only on $n\in\NN$,
such that the function $\cV^n$ in \cref{ET3.2A} satisfies
\begin{equation*}
\widehat\Ag^n_z \cV^n(\Hat{x},s) \,\le\,
c_0(n) - c_1(n) \cV^n(\Hat{x},s) \qquad\forall\,
(\Hat{x},s)\in \sX^n\times\RR_+^m\,,\ \forall\,z\in\cZn(\Hat{x})\,.
\end{equation*}
In particular, under any work-conserving stationary Markov policy,
the process $(\Hat X^n,S^n)$ is positive Harris recurrent,
and $\norm{\Hat{x}}^{}_1$ is integrable under its invariant probability distribution.
\end{theorem}

\proof{Proof.}
The proof mimics that of \cref{T3.1} also using \cref{R2.6}.
The important difference here is that, if we let
$\Hat\varphi^n(x_i)\df \varphi^n\bigl(\Hat{x}_i^n(x_i)\bigr)$, and
\begin{equation*}
\Hat\phi^n_i(x,s) \,\df\, \bigl(1-\zeta_i^n(s_i)\bigr)
\dd\Hat\varphi^n(x_i;1)\,,
\end{equation*}
then, following the steps in \cref{PL3.1B}, we obtain
\begin{equation}\label{PT3.2A}
\begin{aligned}
\Ag^n_z \Hat\phi^n_i(x,s) &\,=\,
\bigl(\lambda_i^n-r_i^n(s_i)\bigr)\dd\Hat\varphi^n(x_i;1)\\[5pt]
&\mspace{100mu}
 - \bigl(\mu^n_i z_i +\gamma^n_i q_i(x,z)\bigr) \bigl(1-\zeta_i^n(s_i)\bigr)
\bigl(\dd\phi^n(x_i-1;1)-\dd\phi^n(x_i;1)\bigr)\,.
\end{aligned}
\end{equation}
As a result, the terms corresponding to the second line of \cref{PL3.1B},
for which the assumption that the hazard rate functions are bounded
was invoked, are not present in \cref{PT3.2A}.
The rest of the proof is the same.
\Halmos\endproof

Without assuming that the abandonment rates are positive, but with
$\varrho^n>0$,
we obtain uniform stability, that is, tightness of the invariant distributions.
To establish this, we scale the Lyapunov function in \cref{ET3.2A}, with
a parameter $\veo>0$.
More precisely, we define
\begin{equation}\label{ET3.3A}
\cV^n_\veo(\Hat{x},s)\,\df
\sum_{i\in\cI} \bigl(1-\zeta_i^n(s_i)\bigr)
\dd\varphi^n_\veo(\Hat{x}_i;n^{-\nicefrac{1}{2}})
+ \sum_{i\in\cI} \frac{\varphi^n_\veo(\Hat{x}_i)}{\mu_i^n}\,,
\end{equation}
with 
\begin{equation*}
\varphi^n_\veo(y)\,\df\, \theta^n \veo\psi(-y) + \psi_\veo(y)\,, \quad y\in\RR\,.
\end{equation*}
The parameter $\theta^n$ depends on certain bounds which we review next.
First, as we have seen in \cref{E-hC1}, there is a constant
$\Check{C}_1$ such that
\begin{equation*}
\babs{\dd\varphi^n_\veo(x\pm n^{-\nicefrac{1}{2}}e_j;\pm n^{-\nicefrac{1}{2}}e_i)
-\dd\varphi^n_\veo(x;\pm n^{-\nicefrac{1}{2}}e_i)}
\,\le\, \frac{1}{n} \Check{C}_1 \veo(\veo+\theta)
\qquad\forall\,i,j\in\cI\,.
\end{equation*}
Let also $\Check{C}_0^n$ be a bound for $\norm{\max_i\,\zeta_i^n}_\infty$.
With $\widehat{C}^n_2$ the constant in \cref{PT3.1B}, we define
\begin{equation*}
\Bar{C}_0^n\,\df\,m^2 \Check{C}_0^n \Check{C}_1 \,,
\quad\text{and\ \ } \Bar{C}_1^n \,\df\,
\Check{C}_1\Biggl(m^2\Check{C}_0^n \mu_i^n
+ \sum_{i\in\cI} \frac{\lambda_i^n}{n}\Biggr)\,.
\end{equation*}
Let $\theta^n$ be equal to the right-hand side of \cref{PT3.1I}
after we replace $\widehat{C}_1$, $\widetilde{C}^n_1$, and $\widetilde{C}^n_2$
with $\Check{C}_1$, $\Bar{C}^n_1$, and $\Bar{C}^n_2$, respectively.

\begin{theorem}\label{T3.3}
Grant \cref{A3.1}, and, in addition,
assume that hazard rate functions $\{h_i\}_{i\in\cI}$ are locally bounded.
Suppose that $\varrho^n>0$. 
Then there exist a cube $K$ and a constant $C$ depending on $\veo$, $\varrho^n$,
and $\theta^n$, defined above,
such that the function $\cV^n_\veo$ in \cref{ET3.3A} satisfies
\begin{equation*}
\widehat\Ag^n_z \cV^n_\veo(\Hat{x},s) \,\le\,
\veo C\Ind_{K}(\Hat{x}) - \veo \frac{\varrho^n}{3m} \qquad\forall\,
(\Hat{x},s)\in \sX^n\times\RR_+^m\,,\ \forall\,z\in\cZn(\Hat{x})\,,
\end{equation*}
and for all $\veo\le\theta^n$.
In particular, under any work-conserving stationary Markov policy,
the process $(\Hat X^n,S^n)$ is positive Harris recurrent.
\end{theorem}

\proof{Proof.}
We follow the proofs of \cref{L3.1,T3.1} to obtain
the analogous inequality to \cref{PT3.1H}.
The result then follows by applying the drift inequality
in \cref{EL2.1C}, and using the definition
of $\theta^n$.
\Halmos\endproof

\subsection{Results with Poisson arrivals.}\label{S3.3}

In this subsection, we specialize the results
to a sequence of queueing models with 
Poisson arrivals with rates $\lambda^n_i$, $i\in\cI$. 
Here, under a stationary Markov policy,
the process $\process{X^n}$ is Markov with extended generator
\begin{equation}\label{E-sA}
\sA^n_z f(x) \,\df\, 
\sum_{i\in\cI} \lambda^n_i \bigl(f(x+e_i) - f(x)\bigr)
+ \sum_{i\in\cI}\bigl(\mu^n_i z_i +\gamma^n_i q_i(x,z)\bigr)
\bigl(f(x-e_i) - f(x)\bigr)\,.
\end{equation}
Define $\widehat\sA^n_z$ analogously to \cref{E-hatAg}.
Mimicking the proof of \cref{T3.1}, we deduce the following,
which we state without proof.

\begin{corollary}\label{C3.1}
Assume that the arrival processes are Poisson. 
Suppose $\varrho^n>0$.
Then for some $\theta=\theta(n)>0$, there exist
positive constants $\Hat\veo_0(n)$ and $C_0^n(\veo)$, such that
the function $V$ in \cref{EL2.1B} satisfies
\begin{equation*}
\widehat\sA^n_z V(\Hat{x}) \,\le\,
C_0^n(\veo) - \veo\frac{\varrho^n}{2m} V(\Hat{x}) \qquad\forall\,
\Hat{x}\in \sX^n\,,\ \forall\,z\in\cZn(\Hat{x})\,,
\end{equation*}
and for all $\veo\in(0,\Hat\veo_0(n))$.
In particular, under any work-conserving stationary Markov policy,
the process $\process{\Hat X^n}$ is exponentially ergodic,
and $V(\Hat{x})$ is integrable under its invariant probability measure.
\end{corollary}

\begin{remark}
Let $\fZn$ denote the class of work-conserving stationary Markov policies
for the process $\Hat{X}^{n}(t)$.
Suppose $\varrho>0$, and let $P_t^{n,z}$ and $\uppi^n_z$
denote the transition probability and the stationary distribution, respectively, of
$\Hat{X}^{n}(t)$ under a policy $z\in\fZn$.
Then, \cref{C3.1} implies
that there exist positive constants $\gamma$ and $C_\gamma$ not depending on $n$ or $z$,
such that
\begin{equation}\label{E5a}
\bnorm{P^{n,z}_t(\Hat{x},\cdot\,)-\uppi^n_z(\cdot)\,}_{V}\,\le\,
C_\gamma V(\Hat{x})\, \E^{-\gamma t}\,,\qquad \forall\,\Hat{x}\in\sX^n\,,
\ \forall\,t\ge0\,.
\end{equation}
Also
\begin{equation*}
\sup_{n\in\NN}\;\sup_{z\in\fZn}\; \int_{\sX^n} V(\Hat{x})\,\uppi^n_z(\D\Hat{x})
\,<\,\infty\,.
\end{equation*}

Note that, if $\nu^n$ denotes the distribution of $\Hat{X}^{n}(0)$, then
\cref{E5a} implies that
\begin{equation}\label{E5b}
\bnorm{P^{n,z}_t(\nu^n,\cdot\,)-\uppi^n_z(\cdot)\,}_{V}\,\le\,
C_\gamma \nu^n(V)\, \E^{-\gamma t}\qquad\forall\,t\ge0\,,
\end{equation}
where $P^{n,z}_t(\nu^n,\cdot\,)\df
\int_{\sX^n}\nu^n(\D\Hat{x}) P^{n,z}_t(\Hat{x},\cdot\,)$
and $\nu^n(V)\df \int_{\sX^n} V(\Hat{x})\nu^n(\D\Hat{x})$.
In particular,  if $\Hat{X}^{n}(0)$ is such that
$\sup_{n\in\NN}\,\nu^n(V)<\infty$, then the convergence in
\cref{E5b} is uniform over  $z\in\fZn$ and $n\in\NN$.
\end{remark}

We also wish to remark that, provided that the jobs do not abandon the queues, that is, 
$\varGamma=0$, the hypothesis $\varrho^n>0$ is sharp.
In fact, there is a dichotomy.
As shown in \cref{C3.1}, if $\varrho^n>0$,
then $\process{X^n}$ is uniformly exponentially ergodic.
Following for example the proof in \cite[Theorem~3.3]{APS18}
one can show that if $\varrho^n<0$ and jobs do not abandon the queues, then
$\process{X^n}$ is transient under any Markov scheduling policy.

As explained in \cite[p.~33]{GS-12}, under positive abandonment in all classes,
the invariant distribution of $\Hat{X}^n$
cannot integrate a function of the form $\E^{\veo \abs{\Hat{x}}^2}$ for $\veo>0$,
even though the invariant probability distribution of the limit diffusion has
this property as seen in \cref{T2.2}.
Note that the technique in the proof of \cref{T3.1} stumbles in \cref{PT3.1D}, since
this bound is no longer valid for the function $\widetilde{V}$ of \cref{T2.2}.

Nevertheless, we have the following improvement of \cref{C3.1},
under positive abandonment in all classes.

\begin{theorem}\label{T3.4}
Assume that the arrival processes are Poisson. 
Suppose $\liminf_{n\to\infty}\,\gamma_i^n>0$ for all $i\in\cI$.
Then there exist positive constants
$\Breve\kappa_0(\eta)$ and $\Breve\kappa_1(\eta)$, such that the function
\begin{equation*}
\Breve{V}^n(\Hat{x}) \,\df\, \exp\bigl(\Phi^*_{\eta,\theta^n}(\Hat{x})\bigr)
\,=\, \exp\bigl(\eta\theta^n\Psi(-\Hat{x})+\eta\Psi(\Hat{x})\bigr)\,,
\end{equation*}
with $\theta^n$ given by \cref{ET3.2B}, satisfies
\begin{equation*}
\widehat\sA^n_z \Breve{V}^n(\Hat{x}) \,\le\,
\Breve\kappa_0(\eta) - \Breve\kappa_1(\eta)\norm{\Hat{x}}^{}_1 \Breve{V}^n(\Hat{x})
\qquad\forall\, (\Hat{x},z)\in\sX^n\times\cZn(\Hat{x})\,,
\end{equation*}
and for all sufficiently large $n$.
In particular, the function $\exp\bigl(\eta\norm{\Hat{x}^n}^{}_1\bigr)$ is integrable
under the stationary distribution of $\process{\Hat{X}^n}$ for all $\eta>0$,
under any work-conserving stationary Markov scheduling policy.
\end{theorem}

\proof{Proof.}
Let $\widehat{\cV}^n(x)\df\Breve{V}\bigl(\Hat{x}^n(x)\bigr)$.
Applying the operator in \cref{E-sA} to $\widehat{\cV}^n$ and using
the analogous bound to \cref{PT3.1D},
\begin{align*} 
\sA^n_z\, \widehat{\cV}^n(x) &\,\le\,
\sum_{i\in\cI} \biggl[\lambda^n_i
\Bigl(\partial_{x_i} \widehat{\cV}^n(x)
+\tfrac{1}{n}\eta(1+\theta)\,\widehat{C}\,\widehat{\cV}^n(x)\Bigr)\nonumber\\
&\mspace{150mu}
 + \bigl(\mu^n_iz_i+\gamma^n_i q_i(x,z)\bigr)\,
\Bigl(-\partial_{x_i} \widehat{\cV}^n(x)
+\tfrac{1}{n}\eta(1+\theta)\,\widehat{C}\,\widehat{\cV}^n(x)\Bigr) \biggr]
\end{align*}
for some constant $\widehat{C}$.
Using \cref{PL3.1D} we write this as
\begin{align}\label{PT3.4A} 
\widehat\sA^n_z\, \Breve{V}(\Hat{x}) &\,\le\,
\sum_{i\in\cI} \biggl(-\frac{\varrho^n \mu^n_i}{d}
 - \mu^n_i\Hat{z}_i - \gamma^n_i \Hat{q}_i(\Hat{x},\Hat{z})\biggr)\,
 \partial_{\Hat{x}_i} \Breve{V}(\Hat{x})\nonumber\\[5pt]
 &\mspace{150mu} +
\eta(1+\theta)\,\widehat{C}\, \sum_{i\in\cI} \biggl(\frac{\lambda^n_i}{n}
+\frac{\mu^n_i}{n}z_i
+\frac{1}{\sqrt n}\gamma^n_i(\Hat{x}_i-\Hat{z}_i)\biggr) \,\Breve{V}(\Hat{x})\,.
\end{align}
Thus, using the drift inequality in \cref{R2.6} to bound
the first term on the right-hand side of \cref{PT3.4A}, and noting
that the coefficient of $\Breve{V}$ on the second term on the right-hand
side is of order $\frac{1}{\sqrt n} \norm{x}^{}_1$,
we establish the result.
\Halmos\endproof

We conclude with the analogous result to \cref{C2.1}.
We need the following notation.
\begin{equation*}
\Hat\cI_1\,\df\,\Bigl\{i\in\cI\colon \limsup_{n\to\infty}\tfrac{\gamma_i^n}{\mu_i^n}<1
\Bigr\}\,.
\end{equation*}

\begin{theorem}\label{T3.5}
Assume that the arrival processes are Poisson. 
Suppose $\liminf_{n\to\infty}\,\varrho^n>0$.
Then the function
$\exp\Bigl(\eta\sum_{i\in\Hat\cI_1} \Hat{x}_i^n\Bigr)$
is integrable under the invariant probability distribution
of $\process{\Hat{X}^n}$ for all $\eta>0$,
and for all sufficiently large $n$.
\end{theorem}

The proof closely mimics that of \cref{T3.1}, and is therefore omitted.

\section*{Acknowledgment.}
This work was supported in part by an Army Research Office grant W911NF-17-1-0019,
in part by NSF grants DMS-1715210, CMMI-1635410, and DMS/CMMI-1715875,
and in part by the Office of Naval Research through grant N00014-16-1-2956.


\end{document}